\documentclass[11pt,twoside,a4paper]{article}
\usepackage[french,english]{babel} 
\usepackage{a4} 
\usepackage[T1]{fontenc} 
\usepackage{epsfig} 
\usepackage{amsmath, amsthm} 
\usepackage{amsfonts,amssymb}
\usepackage{float} 
\usepackage{amsmath,amssymb}
\usepackage{dsfont}
\usepackage[round]{natbib}
\usepackage[usenames]{color}

\newtheorem{theorem}{Theorem}
\newtheorem{main}{Main Lemma}

\newtheorem{lemma}[theorem]{Lemma}

\newtheorem{proposition}[theorem]{Proposition}

\newcounter{as}
\newtheorem{assum}[as]{Assumption}

\DeclareMathOperator*{\argmin}{arg\,min}
\DeclareMathOperator*{\argmax}{arg\,max}

\title{\huge \bf Conditional density estimation in a censored single-index regression model}
\author{OLIVIER BOUAZIZ$^*$ and OLIVIER LOPEZ}
\date{}
\begin{document}
\maketitle 
\noindent Laboratoire de Statistique Théorique et Appliquée,\\ 
	Université Paris VI, \\
  175 rue du Chevaleret \\
  75013 Paris, France\\


\thispagestyle{empty}

\begin{abstract}
Under a single-index regression assumption, we introduce a new semiparametric procedure to estimate a conditional density of a censored response. The regression model can be seen as a generalization of Cox regression model and also as a profitable tool to perform dimension reduction under censoring. This technique extends the results of Delecroix {\it et al.} (2003). 
We derive consistency and asymptotic normality of our estimator of the index parameter by proving its asymptotic equivalence with the (uncomputable) maximum likelihood estimator, using martingales results for counting processes and arguments of empirical processes theory. Furthermore, we provide a new adaptive procedure which allows us both to chose the smoothing parameter involved in our approach and to circumvent the weak performances of Kaplan-Meier estimator (1958) in the right-tail of the distribution. Through a simulation study, we study the behavior of our estimator for small samples.\\


\end{abstract}

\section{Introduction}

A major issue of recent papers dealing with censored regression is to propose alternatives to the popular Cox regression model. This model, also known as multiplicative hazard regression model (see \cite{Cox72}), states some semiparametric model on the conditional hazard function. Estimation in this model is traditionally performed using pseudolikelihood techniques, and the theoretical properties of these procedures are covered by a large number of papers (see e.g. \cite{Flemming91}). However, in some situations, the assumptions of Cox regression model are obviously not satisfied by the data set. In this paper, our aim is to perform estimation in a semiparametric regression model which allows more flexibility than the Cox regression model. This new technique can be seen as a particularly interesting alternative, since it is valid in a larger number of situations than the multiplicative hazard model.

Alternatives to Cox regression model mostly focus on the estimation of a conditional expectation, or of a quantile regression model. Koul {\it et al.} (1981), \cite{Stute99}, Delecroix {\it et al.} consider mean-regression models where the regression function belongs to a parametric family, but with an unknown distribution of the residuals. Parametric quantile regression was studied by Gannoun {\it et al}. On the other hand, \cite{Lu05} and \cite{newLopez} considered a semiparametric single-index regression model. Single-index regression models were initially introduced to circumvent the so-called ``curse of dimensionality'' in nonparametric regression (see. e.g. \cite{Ichimura93}), by assuming that the conditional expectation only depends on an unknown linear combination of the covariates. Another appealing aspect of such kind of models is that they include the Cox regression model as a particular case. The main assumption of this model is that the conditional density only depends on an unknown linear combination of the covariates, while the multiplicative hazard model states a similar assumption on the conditional hazard rate. In this paper we focus on estimation of the parameter in a regression model in which the conditional density of the response satisfies a single-index assumption. We provide asymptotic results for a new M-estimation procedure for the index parameter. This procedure can be seen as a generalization of the method of Delecroix {\it et al.} (2003) to the case of censored regression.




As in the uncensored case, we show that, regarding to the estimation of the parametric part of our model, there is an asymptotic equivalence between our semiparametric approach and a parametric one relying on some prior knowledge on the family of regression functions. For the nonparametric part, we use kernel estimators of conditional densities as in Delecroix {\it et al.} (2003). Since the performance of kernel estimators strongly relies on the choice of the smoothing parameter, we also provide a method to choose this parameter adaptively. Another technical issue in our approach concerns a truncation parameter involved in our procedure. This problem of truncation directly comes from the censored framework, where estimators of the underlying distribution functions sometimes fail to estimate correctly the tail of the distribution. This problem is traditionally circumvented for example by assuming integrability assumptions on the response and censoring distribution, see e.g. \cite{Stute99}. On the other hand, truncation procedure consists of removing the observations which are too large in the estimation of the regression function, see e.g. \cite{Heuchenne07}, or condition (2.2) in \cite{Brunel06} which can be interpreted as such kind of truncation. Until now, the truncation bounds which were used were arbitrary fixed, and usually no method is proposed to discuss a method for choosing this truncation bound in practical situations. Therefore, in the new method we propose, we also provide a data-driven procedure to choose the truncation parameter. In our practical implementations, we used a criterion based on an asymptotic discussion which focuses on the mean-squared error associated with the estimation of the single index parameter. We also suggest some possible adaptations to other type of criterion which are covered by our theoretical results.

In section \ref{secreg}, we introduce our censored regression model and present our estimation procedure. It relies on the Kaplan-Meier estimator (1958) of the distribution function, and on semiparametric estimators of the conditional density. Following the procedure of Delecroix {\it et al.} (2003), we considered kernel based estimators. Our theoretical results are presented in section \ref{theorique}. In section \ref{sec_simul} we report simulation results and analysis on real data. Section \ref{lemme} contains the detailed proof of our Main Lemma which states the asymptotic equivalence of estimating the parameter in the semiparametric and parametric models. All the technicalities are postponed to the section \ref{technical}.


\section{Censored regression model and estimation procedure}
\label{secreg}

\subsection{Notations and general setting} 
\label{secnot}

Let $Y_1, \ldots, Y_n$ be i.i.d. copies of a random response variable $Y
\in \mathbb{R},$ and let $X_1, \ldots,X_n$ be i.i.d. copies of a random
vector of covariates $X\in \mathcal{X},$ where $\mathcal{X}$ is a compact subset of $\mathbb{R}^d$. Introducing $C_1,\ldots,C_n$ i.i.d. replications of the censoring variable $C\in \mathbb{R},$ we consider the following censored regression framework, where the observations are
\[\begin{cases}
Z_i = Y_i\wedge C_i & 1\leq i\leq n,\\
\delta_i \;= \mathds{1}_{\{Y_i\leq C_i\}} & 1\leq i\leq n, \\
X_i  \in {\mathcal{X} \subset \mathbb{R}^d}  & 1\leq i \leq n.
\end{cases}\]
Let us introduce some notations for the distribution functions of the random variables appearing in this model, that is $H(t)= \mathbb{P}(Z\leq t),$ $F_X(t)= \mathbb{P}(X\leq t),$ $F_Y(t)= \mathbb{P}(Y\leq t),$ $F_{X,Y}(x,y)= \mathbb{P}(X\leq x,Y\leq y)$
and $G(t)= \mathbb{P}(C\leq t).$ A major difficulty arising in censored regression models stands in the unavailability of the empirical distribution function to estimate functions $F_Y,$ $F_{X,Y}$ and $G,$ which must be replaced by Kaplan-Meier estimators.

We are interested in estimating $f(y|x),$ where $f(y|x)$ denotes the conditional density of $Y$ given $X=x$ evaluated at point $y.$ If one has no insight on the function $f,$ it becomes necessary to perform nonparametric estimation of the conditional density. In absence of censoring, a classical way to proceed is to use kernel smoothing, see e.g. \cite{Bashtannyk01}. However, the so-called ``curse of dimensionality'' prevents this approach from being of practical interest when the number of covariates is important ($d>3$ in practice). Therefore it becomes relevant to consider semiparametric models which appear to be a good compromise between the parametric (which relies on strong assumptions on the function $f$ which may not hold in practice) and the nonparametric approach (which relies on fewer assumptions). In the following, we will consider the following semiparametric single-index regression model,
\begin{equation}
\label{model} \exists \, \theta_0 \in \Theta \subset\mathbb{R}^ d \; s.a. \; f(y|x)=f_{\theta_0}(y,x'\theta_0),
\end{equation}
where $f_{\theta}(y,u)$ denotes the conditional density of $Y$ given $X'\theta=u$ evaluated at $y.$ For identifiability reasons, we will impose that the first component of $\theta_0$ is one. In comparison to Cox regression model for absolute continuous variables, our model (\ref{model}) is more general, since it only assumes that the law of $Y$ given $X$ depends on an unknown linear combination of the covariates, without imposing additional conditions on the conditional hazard rate.

Model (\ref{model}) has been considered by Delecroix {\it et al}. (2003) in the uncensored case. However, their procedure can not be directly applied in the censored framework since the responses variables are not directly observed. As a consequence, the empirical distribution function is unavailable, and most of the tools used in this context are not at our disposal. A solution consists of using procedures relying on Kaplan-Meier estimators for the distribution function. An important difficulty arising in this type of techniques stands in the poor behavior of Kaplan-Meier estimators in the tail of the distribution. A practical way to prevent us from this kind of drawback is to consider truncated version of the variable $Y.$ In the following, we will consider $A_{\tau}$ a sequence of compacts included in the set $\{t:\tau_1\leq t \leq \tau\},$ for $\tau\leq \tau_0,$ where $\tau_0<\inf\{t:H(t)=1\}.$ Using only the observations in $A_{\tau}$ allow us to avoid the bad behavior of usual Kaplan-Meier estimators in the tail of the distribution. Moreover, this technique of truncation is particularly adapted to our problem of estimating $\theta_0.$ In our framework, this truncation does not lead to any asymptotic bias, since, denoting by $f^{\tau}(\cdot|x)$ the conditional density of $Y$ given $X=x$ and $Y\in A_\tau,$ for any $\tau<\infty,$ we have, under (\ref{model}),
\begin{equation}
\label{model2}
f^{\tau}(y|x)=f^{\tau}_{\theta_0}(y,x'\theta_0),
\end{equation}
where $f^{\tau}_{\theta}(y,u)$ denotes the conditional density of $Y$ given $X'\theta=u$ and $Y\in A_\tau$ evaluated at $y,$
and where the parameter is the same in (\ref{model}) as in (\ref{model2}). In section \ref{tauadapt}, we will discuss a new method allowing to choose $\tau$ from the data in order to improve the performance in estimating $\theta_0$.

%


\subsection{Estimation procedure}

We will extend the idea behind the procedure developed by Delecroix {\it et al.} (2003), adapting it to our censored framework.
First assume that we know the family of functions $f^{\tau}_{\theta}.$ This approach is a modification of the maximum likelihood estimation procedure. Define, for any function $J\geq 0,$
\begin{align*}
L^{\tau}(\theta,J) &= E\left[\log f^{\tau}_{\theta}(Y_i,\theta'X_i)J(X_i)\mathds{1}_{Y_i\in A_{\tau}}\right]=\int \log f^{\tau}_{\theta}(y,\theta'x)J(x)\mathds{1}_{y\in A_{\tau}}dF_{X,Y}(x,y).
\end{align*}
Here, $J$ is a positive trimming function which will be defined later in order to avoid denominators problems in the nonparametric part of the model, see section \ref{J}. From (\ref{model2}), $\theta_0$ maximizes $L^{\tau}(\theta,J)$ for any $\tau<\infty,$ this maximum being unique under some additional conditions on the regression model and $J.$ Since, in our framework, $F_{X,Y}$ and $f^{\tau}_{\theta}$ are unknown, it is natural to estimate them in order to produce an empirical version of $L^{\tau}(\theta,J).$

\subsubsection{Estimation of $F_{X,Y}$} In the case where there is no censoring (as in Delecroix {\it et al.} (2003)), $F_{X,Y}$ can be estimated by the empirical distribution function. 
In our censoring framework, the empirical distribution function of $(X,Y)$ is unavailable, since it relies on the true $Y_i'$s which are not observed. A convenient way to proceed consists of replacing it by some Kaplan-Meier estimator such as the one proposed by \cite{Stute93}. Let us define the Kaplan-Meier estimator (\cite{Kaplan58}) of $F_Y,$
\begin{align*}
\hat{F}_Y(y) &= 1-\prod_{i:Z_i\leq t}\left(1-\frac{1}{\sum_{j=1}^n\mathds{1}_{Z_j\geq Z_i}}\right)^{\delta_i}
\\
&=\sum_{i=1}^n \delta_i W_{in}\mathds{1}_{Z_i\leq y},
\end{align*}
where $W_{in}$ denotes the ``jump'' of Kaplan-Meier estimator at observation $i$ (see \cite{Stute93}). To estimate $F_{X,Y},$ Stute proposes to use
\[\hat{F}(x,y)=\sum_{i=1}^n \delta_iW_{in}\mathds{1}_{Z_i\leq y,X_i\leq x}.\]
Let us also define the following (uncomputable) estimator of the distribution function,
\[\tilde{F}(x,y)=\sum_{i=1}^n \delta_iW_{i}^*\mathds{1}_{Z_i\leq y,X_i\leq x},\]
where $W_i^*=n^{-1}[1-G(Z_i-)]^{-1}.$ The link between $\hat{F}$ and $\tilde{F}$ comes from the fact that, in the case where $\mathbb{P}(Y=C)=0,$ \begin{equation}\label{explicit_jump}
W_{in}=n^{-1}[1-\hat{G}(Z_i-)]^{-1},
\end{equation} where $\hat{G}$ denotes the Kaplan-Meier estimator of $G$ (see \cite{Satten01}). Asymptotic properties of $\hat{F}$ can be deduced from studying the difference with the simplest but uncomputable estimator $\tilde{F}.$

If we know the family of regression functions $f^{\tau}_{\theta},$ it is possible to compute the empirical version of $L^{\tau}(\theta,J)$ using $\hat{F},$ that is
\begin{align*}
L_n^{\tau}(\theta,f^{\tau},J) &= \int \log f^{\tau}_{\theta}(y,\theta'x)J(x)\mathds{1}_{y\in A_\tau}d\hat{F}(x,y) \\
&= \sum_{i=1}^n \delta_iW_{in}\log f^{\tau}_{\theta}(Z_i,\theta'X_i)J(X_i)\mathds{1}_{Z_i\in A_{\tau}}.
\end{align*} In the case $J\equiv 1,$ the estimator of $\theta_0$ obtained by maximizing $L^{\tau}_n$ would turn out to be an extension of the maximum likelihood estimator of $\theta_0,$ used in presence of censoring.

\subsection{Estimation of $f^{\tau}_{\theta}$}

In our regression model (\ref{model2}), the family $\{f_{\theta}^{\tau},\theta\in \Theta\}$ is actually unknown. As in Delecroix {\it et al.} (2003), we propose to use nonparametric kernel smoothing to estimate $f_{\theta}^{\tau}.$ Introducing a kernel function $K$ and a sequence of bandwidths $h,$ define
\begin{align}
\label{kernel_estimator} \hat{f}_{\theta}^{h,\tau}(z,\theta'x) &= \frac{\int
K_h(\theta'x-\theta'u)K_h(z-y)\mathds{1}_{y\in A_{\tau}}d\hat{F}(u,y)}{\int
K_h(\theta'x-\theta'u)\mathds{1}_{y\in A_{\tau}}d\hat{F}(u,y)},
\end{align}
where $K_h(\cdot)=h^{-1}K(\cdot/h).$ Also define $f^{*h,\tau}$ the kernel estimator based on function $\tilde{F},$ that is
\begin{align*}
f^{*h,\tau}_{\theta}(z,\theta'x) &= \frac{\int
K_h(\theta'x-\theta'u)K_h(z-y)\mathds{1}_{y\in A_{\tau}}d\tilde{F}(u,y)}{\int
K_h(\theta'x-\theta'u) \mathds{1}_{y\in A_{\tau}}d\tilde{F}(u,y)}.
\end{align*}
$f^{*h,\tau}$ will play an important role in studying the asymptotic behavior of $\hat{f}^{h,\tau}.$ Indeed, $f^{*h,\tau}$ is theoretically more easy to handle with, since it relies on sums of i.i.d. quantities, which is not the case for $\hat{F}.$ Since $f^{*h,\tau}$ can be studied by standard kernel arguments, the most important difficulty will arise from studying the difference between $\hat{f}^{h,\tau}$ and $f^{*h,\tau}.$

In the following, we will impose the conditions below on the kernel function.
\begin{assum}
\label{noyau}Assume that
\begin{itemize}
\item [(A1)] $K$ is a twice differentiable and four order kernel with derivatives of order 0, 1 and 2 of bounded variation. Its support is contained in $[-1/2,1/2]$ and $\int_{\mathbb{R}}K(s)ds =1$,
\item [(A2)] $\|K\|_{\infty}:=\sup_{x\in \mathbb{R}}|K(x)|<\infty$,
\item [(A3)] $\mathcal{K}:=\{K\big((x-\cdot)/h\big):h>0, x\in \mathbb{R}^d\}$ is a pointwise measurable class of functions,
\item [(A4)] $h\in \mathcal{H}_n\subset [an^{-\alpha},bn^{-\alpha}]$ with $a,b\in \mathbb{R},$ $1/8<\alpha<1/6$ and where $\mathcal{H}_n$ is of cardinality $k_n$ satisfying $k_nn^{-4\alpha}\rightarrow 0.$
\end{itemize}
\end{assum}

\subsection{The trimming function $J$}
\label{J}

The reason behind introducing function $J$ has to be connected with the need to prevent us from denominators close to zero in the definition (\ref{kernel_estimator}). Ideally, we would need to use the following trimming function,
\begin{equation}
\label{J0} J_0(x,c)=\tilde{J}(f_{\theta_0'X},\theta_0'x,c),
\end{equation}
where $c$ is a strictly positive constant, $f_{\theta_0'X}$ denotes the density of $\theta_0'X$ and $\tilde{J}(g,u,c)=\mathds{1}_{g(u)>c}.$
Unfortunately, this function relies on the knowledge of parameter $\theta_0$ and $f_{\theta_0'X}.$ Therefore, we will have to proceed in two steps, that is first obtain a preliminary consistent estimator of $\theta_0,$ and then use it to estimate the trimming function $J_0$ which will be needed to achieve asymptotic normality of our estimators of $\theta_0.$

We will assume that we know some set $B$ on which $\inf\{f_{\theta'X}(\theta'x): x\in B, \theta \in \Theta\}>c,$ where $c$ is a strictly positive constant. In a preliminary step, we can use this set $B$ to compute the preliminary trimming $J_B(x)=\mathds{1}_{x\in B}.$ Using this trimming function, and a deterministic sequence of bandwidth $h_0$ satisfying $(A4)$ in Assumption \ref{noyau}, we define a preliminary estimator $\theta_n$ of $\theta_0,$
\begin{equation}
\label{preliminary} \theta_n =\argmin_{\theta\in \Theta}L^{\tau}_n(\theta,\hat{f}^{h_0,\tau},J_B).
\end{equation}
Let us stress the fact that $B$ is assumed to be known by the statistician. This is a classical assumption in
single-index regression (see Delecroix {\it et al.}, (2006)). However, in practice, the procedure does not seem
very sensitive to the choice of $B.$ The bandwidth $h_0$ we consider in the preliminary step can be any sequence decreasing
to zero slower than $n^{-1/2}.$ Adaptive choice of $h_0$ could be considered (using, for instance, the same choice as in the final estimation step, see below). However, since we will only need $\theta_n$ to be a preliminary consistent estimator, and the final estimator will not be very
sensitive to an adaptive choice of $h_0$ while computing $\theta_n,$ we do not consider this case in the following.

With, at hand, this preliminary estimator $\theta_n,$ we can compute an estimated version of $J_0$ which will happen to be equivalent to $J_0$ (see Delecroix {\it et al.} (2006) page 738), that is
\begin{equation}
\label{Jhat} \hat{J}_0(x,c)=\tilde{J}(\hat{f}^{h_0,\tau}_{\theta_n'X},\theta_n'x,c).
\end{equation}

For each sequence of bandwidths satisfying (A4) in Assumption \ref{noyau}, and for each truncation bound $\tau,$ we can define an estimator of $\theta_0$
\begin{equation}
\label{thetahat} \hat{\theta}^{\tau}(h)=\argmax_{\theta \in \Theta_n} L^{\tau}_n(\theta,\hat{f}^{h,\tau},\hat{J}_0),
\end{equation}
where $\Theta_n$ is a shrinking sequence of neighborhoods accordingly to the preliminary estimation.
However, as for any smoothing approach, the performance of this procedure strongly depends on the bandwidth sequence. Therefore it becomes particularly relevant to provide an approach which automatically selects the most adapted bandwidth according to the data. Then, the new question arising from the censored framework comes from the adaptive choice of the truncation parameter $\tau.$

\subsection{Adaptive choice of the bandwidth}
\label{hadapt}

Our procedure consists of choosing from the data, for each $\theta,$ a bandwidth
which is adapted to the computation of $f^{\tau}_{\theta}(z,u).$ For
this, we use an adaptation of the cross-validation technique of \cite{Fan04}, that is
$$\hat{h}^{\tau}(\theta)=\argmin_{h \in \mathcal{H}_n}\sum_{i=1}^n W_{in}\mathds{1}_{Z_i\in A_{\tau}}\left\{\int_{A_{\tau}} \hat{f}^{h,\tau}_{\theta}(z,\theta'X_i)^2dz-2\hat{f}^{h,\tau}_{\theta}(Z_i,\theta'X_i)\right\}.$$
This criterion is (up to a quantity which does not depend on $h$) an empirical version of the ISE criterion defined in (3.3) in \cite{Fan04} (in a censored framework),
that is
$\int_{A_{\tau}}\int\{\hat{f}^{h,\tau}_{\theta}(z,\theta'x)-f^{\tau}_{\theta}(z,\theta'x)\}^2f_{\theta'X}(\theta'x)dxdz.$

The estimator of $\theta_0$ with an adaptive bandwidth is now defined
as
\begin{equation}
\label{thetahat2} \hat{\theta}^{\tau}=\argmax_{\theta \in
\Theta_n} L_n^{\tau}(\theta,\hat{f}^{\hat{h},\tau}, \hat{J}_0).
\end{equation}
In the above notation, $\hat{h}$ depends on $\theta$ and $\tau,$ which was
not emphasized to shorten the notation.

\subsection{Adaptive choice of $\tau$}
\label{tauadapt}

As we already mentioned, the Kaplan-Meier estimator does not behave well in the tail of the distribution. For example, if some moment conditions are not satisfied, it is not even $n^{1/2}$-consistent. Moreover, even in the case where an appropriate moment condition holds, it may happen (at least for a finite sample size) that the weights corresponding to the large observations are too important and considerably influence the estimation procedure. For this reason, we introduced a truncation by a bound $\tau.$ However, a large number of existing procedure which also rely on such kind of truncation do not consider the problem of choosing $\tau$ from the data. We propose to select $\tau$ from the data in the following way. Suppose that we have a consistent estimator of the asymptotic mean squared error,
$$E^2(\tau)=\limsup_n E\left[\|\hat{\theta}^{\tau}(\hat{h}^{\tau})-\theta_0\|^2\right],$$ say $\hat{E}^2(\tau)$ satisfying
\begin{equation}
\label{consist_var} \sup_{\tau_1\leq\tau\leq \tau_0}|\hat{E}^2(\tau)-E^2(\tau)|\rightarrow 0, \text{ in probability}.
\end{equation}
Such an estimator will be proposed in section \ref{sec_simul}.
Using this empirical estimator, we propose to choose $\tau$ in the following way, that is
$$\hat{\tau}=\argmin_{\tau_1\leq \tau \leq \tau_0}\hat{E}^2(\tau).$$

Our final estimator of $\theta_0$ is based on an adaptive bandwidth and an adaptive choice of truncation parameter $\tau,$ that is
\[\hat{\theta}=\hat{\theta}^{\hat{\tau}}.\]

As we already said, truncating the data does not introduce additional bias in the estimation of $\theta_0$. On the other hand, removing too many data points could strongly increase the variance and removing some of the largest data points will decrease it. Then, our selection procedure $\hat{\tau}$ is based on estimating the variance of $\hat{\theta}$ and consists of taking from the data the truncation parameter $\tau$ that seems to be the best compromise between these two aspects.

\section{Asymptotic results}{\label{theorique}}
\subsection{Consistency}

The assumptions needed for consistency can basically be split into three categories, that is identifiability assumptions, assumptions on the regression model (\ref{model2}) itself and finally assumptions on the censoring model.

\textbf{Identifiability assumption and assumption on the regression model.}

\begin{assum}
\label{identif} Assume that for all $\tau_1\leq \tau\leq \tau_0$ and all $\theta\in \Theta-\{\theta_0\},$
$$L_{\tau}(\theta_0,J_B)-L_{\tau}(\theta,J_B)>0.$$
\end{assum}

\begin{assum}
\label{lgnu} Assume that for $\theta_1, \theta_2 \in \Theta$, for a bounded function $\Phi(X)$ and for some $\gamma>0$, we have
$$\sup_{\tau}\|f^{\tau}_{\theta_1}(y,\theta_1'x)-f^{\tau}_{\theta_2}(y,\theta_2'x)\|_{\infty} \leq \|\theta_1-\theta_2\|^{\gamma}\Phi(X).$$
\end{assum}



\textbf{Assumptions on the censoring model.}

\begin{assum}
$\mathbb{P}(Y=C)=0.$
\end{assum}

This classical assumption in a censored framework avoids problems caused by the lack of symmetry between $C$ and $Y$ in the case where there are ties.

\begin{assum}
\label{ident} Identifiability assumption : we assume that
\begin{itemize}
\item[-] $Y$ and $C$ are independent.
\item[-] $\mathbb{P}(Y\leq C|X,Y)=\mathbb{P}(Y\leq C|Y).$
\end{itemize}
\end{assum}
This last assumption was initially introduced by \cite{Stute93}. An
important particular case in which assumption \ref{ident} holds is when
$C$ is independent from $(X,Y).$ However, assumption \ref{ident} is a more
general and widely accepted assumption, which allows the
censoring variables to depend on the covariates.


\begin{theorem}
\label{consistency} Under Assumptions \ref{identif} to \ref{ident},
\begin{equation}
\label{consist1} \sup_{\theta\in \Theta,{\tau_1\leq \tau \leq \tau_0} }|L_n^{\tau}(\theta,\hat{f}^{h_0,\tau},J_B)-L^{\tau}(\theta,J_B)|=o_P(1),
\end{equation}
and consequently,
$$\theta_n \rightarrow_{\mathbb{P}}\theta_0.$$
\end{theorem}

\begin{proof}
To show (\ref{consist1}), we will proceed in two steps. First we consider $L_n^{\tau}(\theta,f^{\tau},J_B)-L^{\tau}(\theta,J_B)$ (parametric problem), and then $L_n^{\tau}(\theta,\hat{f}^{h_0,\tau},J_B)-L_n^{\tau}(\theta,f^{\tau},J_B).$

\textbf{Step 1.} From Assumption \ref{lgnu}, the family $\{\log(f_{\theta}^{{\tau}}(\cdot,\theta'\cdot )),\theta \in \Theta,\,{\tau_1\leq \tau \leq \tau_0}\}$ is seen to be $P-$ Glivenko-Cantelli. Using an uniform version of \cite{Stute93} leads to $\sup_{\theta}|L_n^{\tau}(\theta,f^{\tau},J_B)-L^{\tau}(\theta,J_{B})|\rightarrow_{\mathbb{P}} 0.$

\textbf{Step 2.} We have, on the set $\Theta'B,$
$$|\log \hat{f}_{\theta}^{h_0,\tau}(y,u)-\log f_{\theta}^{{\tau}}(y,u)|\leq c^{-1}[\hat{f}_{\theta}^{h_0,\tau}(y,u)-f_{\theta}^{{\tau}}(y,u)].$$
Hence,
\begin{align*}
\sup_{\theta,{\tau}}|L_n^{\tau}(\theta,\hat{f}^{h_0,\tau},J_{B})-L_n(\theta,f^{\tau},J_{B})| & \leq  c^{-1}\sup_{\theta,y,u,\tau}|\hat{f}_{\theta}^{h_0,\tau}(y,u)-f_{\theta}^{{\tau}}(y,u)|\mathds{1}_{u\in \Theta'B,y\leq\tau}\int d\hat{F}{(x,y)}
\\
& \leq c^{-1}\sup_{\theta,y,u,\tau}|\hat{f}_{\theta}^{h_0,\tau}(y,u)-f_{\theta}^{{\tau}}(y,u)|\mathds{1}_{u\in \Theta'B,y\leq\tau}.
\end{align*}
Using the uniform convergence of $\hat{f}_{\theta}^{h,\tau}$ (see Proposition \ref{convergence_rate} and Lemma \ref{leplusdur}), deduce that\\ $\sup_{\theta,\tau}|L^{\tau}_n(\theta,\hat{f}^{h_0,\tau},{J_B})-L^{\tau}_n(\theta,f^{\tau},{J_B})|\rightarrow_{\mathbb{P}} 0.$
\end{proof}

\subsection{Asymptotic normality}
\label{sec_normal}

To obtain the asymptotic normality of our estimator, we need to add some regularity assumptions on the regression model.

\begin{assum}
\label{regularity} Denote by $\nabla_{\theta}f_{\theta}^{\tau}(y,x)$ (resp. $\nabla^2_{\theta}f_{\theta}^{\tau}(y,x)$) the vector of partial derivatives  (resp. the matrix of second derivatives with respect to $\theta$) of $f_{\theta}^{\tau}$ with respect to $\theta$ and computed at point $(\theta,x,y).$
Assume that for $\theta_1, \theta_2 \in \Theta$, for a bounded function $\Phi(X)$ and for some $\gamma>0$, we have
$${\sup_\tau} \|\nabla^2_{\theta}f^{\tau}_{\theta_1}(y,x)-\nabla^2_{\theta}f^{\tau}_{\theta_2}(y,x)\|_{\infty} \leq \|\theta_1-\theta_2\|^{\gamma}\Phi(X).$$
\end{assum}

\begin{assum}
\label{donsker1} Using the notation of \cite{Vaart96} in section 2.7, define
\begin{align*}
\mathcal{H}_1 &=
\mathcal{C}^{1+\delta}(\theta_0'\mathcal{X}\times
A_{\tau},M), \\
\mathcal{H}_2 &=
x\mathcal{C}^{1+\delta}(\theta_0'\mathcal{X}\times
A_{\tau},M)+\mathcal{C}^{1+\delta}(\theta_0'\mathcal{X}\times A_{\tau},M)
\end{align*}
Assume that $f^{\tau}_{\theta_0}(\cdot,\cdot)\in \mathcal{H}_1$ (as a function of $\theta_0'x$ and $y$) and $\nabla_{\theta}f^{\tau}_{\theta_0}(\cdot,\cdot)\in \mathcal{H}_2.$
\end{assum}


If the family of functions $f^{\tau}$ {was} known (parametric problem), the asymptotic normality of $\hat{\theta}$ could be deduced from elementary results on Kaplan-Meier integrals (see section \ref{technical} for some brief review of these results), as in \cite{Stute99} or in Delecroix {\it et al.} (2008). Using this kind of results, we can derive the following Lemma (see section \ref{technical} for the proof) which is sufficient to obtain the asymptotic law of $\hat{\theta}$ in the parametric case, from Theorem 1 and 2 of \cite{Sherman94}.

\begin{lemma}
\label{parametric} Under Assumptions \ref{regularity} and \ref{donsker1}, we have the following representations:
\begin{enumerate}
\item on $o_P(1)-$neighborhoods of $\theta_0,$
$$L_n^{\tau}(\theta,f^{ \tau},J_0)=L^{\tau}(\theta,J_0)+(\theta-\theta_0)'T_{1n}(\theta)+(\theta-\theta_0)'T_{2n}(\theta)(\theta-\theta_0)
+T_{3n}(\theta)+T_{4n}(\theta_0),$$
with $\sup_{\theta,{ \tau}}|T_{1n}|=O_P(n^{-1/2}),$ $\sup_{\theta,{ \tau}}|T_{2n}|=o_P(1),$ $\sup_{\theta,{ \tau}}|T_{3n}|=O_P(n^{-1})$ and $T_{4n}(\theta_0)=L^{\tau}_n(\theta_0,f^{ \tau},J_0).$
\item on $O_P(n^{-1/2})-$ neighborhoods of $\theta_0,$
$$L_n^{\tau}(\theta,f^{ \tau},J_0)=n^{-1/2}(\theta-\theta_0)'W_{n,\tau}-\frac{1}{2}(\theta-\theta_0)'V_{\tau}(\theta-\theta_0)+T_{4n}(\theta_0)
+T_{5n}(\theta),$$
\end{enumerate}
with $\sup_{\theta,{ \tau}}|T_{5n}|=o_P(n^{-1}),$
and defining $f_1(x,y)=f_{\theta_0}^{{ \tau}^{-1}}(y,\theta_0'x)J_0(x,c)\nabla_{\theta}f_{\theta_0}^{{ \tau}}(y,x),$
\begin{align*}
W_{n,\tau} &= \frac{1}{n^{1/2}}\sum_{i=1}^n \psi(Z_i,\delta_i,X_i;f_1\mathds{1}_{A_{\tau}}), \\
V_{\tau} &= E\left[f_{\theta_0}^{{ \tau}^{-2}}(Y,\theta_0'X)J_0(X,c)\nabla_{\theta}f_{\theta_0}^{{ \tau}}(Y,X)\nabla_{\theta}f_{\theta_0}^{{ \tau}}(Y,X)'
\mathds{1}_{Y\in A_\tau}\right],
\end{align*}
where $\psi$ is defined in Theorem \ref{km}.
\end{lemma}

In the following Theorem, we show that the semiparametric estimator proposed in section \ref{secreg} has the same asymptotic law as in the fully parametric case.

\begin{theorem}
\label{maintheorem}
Define $\tau^*=\argmin_{\tau} E^{2}(\tau).$
Under Assumptions \ref{noyau} to \ref{donsker1}, we have the following asymptotic i.i.d. representation,
\begin{equation}
\label{iid}
\hat{\theta}-\theta_0= -\frac{1}{n^{1/2}}V_{\tau^*}^{-1}W_{n,\tau^*}+o_P(n^{-1/2}),
\end{equation}
where $V_{\tau}$ and $W_{n,\tau}$ are defined in Lemma \ref{parametric}.
As a consequence,
$$n^{1/2}(\hat{\theta}-\theta_0)\Longrightarrow \mathcal{N}(0,\Sigma_{\tau^*})$$
where $\Sigma_{\tau^*}=V_{\tau^*}^{-1}\Delta_{\tau^*}(f_1) V_{\tau^*}^{-1},$ $\Delta_{\tau^*}(f_1)=Var\left(\psi(Z,\delta,X;f_1\mathds{1}_{A_{\tau^*}})\right)$ and $f_1$ is defined in Lemma \ref{parametric}.
\end{theorem}

This Theorem is a consequence of the Main Lemma 
below. This result shows that, asymptotically speaking, maximizing $L^{\tau}_n(\theta,\hat{f}^{h,\tau},J)$ is equivalent to maximizing $L^{\tau}_n(\theta,f^{ \tau},J).$

\begin{main}
\label{mainlemma} Under Assumptions \ref{noyau} to \ref{donsker1},
$$L^{\tau}_n(\theta,\hat{f}^{h,\tau},\hat{J}_0)=L^{\tau}_n(\theta,f^{ \tau},J_0)+(\theta-\theta_0)'R_{1n}(\theta,h,\tau)+(\theta-\theta_0)'R_{2n}(\theta,h,\tau)(\theta-\theta_0)+\tilde{L}^{\tau}_n(\theta_0),$$
where \begin{align*}
\sup_{\theta \in \Theta_n,h\in \mathcal{H}_n,\tau_1\leq \tau\leq \tau_0} R_{1n}(\theta,h,\tau) &= o_P(n^{-1/2}), \\
\sup_{\theta \in \Theta_n,h\in \mathcal{H}_n,\tau_1\leq \tau\leq \tau_0} R_{2n}(\theta,h,\tau) &= o_P(1).
\end{align*}
and
\[\tilde{L}^{\tau}_n(\theta_0)=A^{\tau}_{1n}(\theta_0,\hat{f}^{h,\tau})-B^{\tau}_{4n}(\theta_0,\hat{f}^{h,\tau})\]
where $A^{\tau}_{1n}(\theta_0,\hat{f}^{h,\tau})$ and $B^{\tau}_{4n}(\theta_0,\hat{f}^{h,\tau})$ are defined in the proof of this Lemma.
\end{main}

In view of Theorem 1 and 2 of \cite{Sherman94}, this result will allow us to obtain the rate of convergence of our estimators, and then the asymptotic law is the same law as the asymptotic law in the parametric problem.

\begin{proof}[Proof of Theorem \ref{maintheorem}]
Define
\begin{align*}
\Gamma_{0n}(\theta,\tau,h) &= L^{\tau}_n(\theta,\hat{f}^{h,\tau},\hat{J}_0), \\
\Gamma_{1n}(\theta,\tau) &= L^{\tau}_n(\theta,\hat{f}^{\hat{h},\tau},\hat{J}_0), \\
\Gamma_{2n}(\theta) &= L^{\hat{\tau}}_n(\theta,\hat{f}^{\hat{h},\hat{\tau}},\hat{J}_0).
\end{align*}
We now apply Theorem 1 and 2 in Sherman (1994) to $\Gamma_{in},$ for $i=0,1,2.$ From our Main Lemma and Lemma \ref{parametric}, we deduce, that the representation (11) in Theorem 2 of \cite{Sherman94} holds for $i=0,1,2,$ on $O_P(n^{-1/2})-$ neighborhoods of $\theta_0,$ with $W_n$ and $V$ defined in Lemma \ref{parametric}.
The asymptotic representation (\ref{iid}) is a by-product of the proof of Theorem 2 in \cite{Sherman94} and of the i.i.d. representations of Kaplan-Meier integrals (see Theorem \ref{km}).
\end{proof}

\section{Simulation study and real data analysis}
\label{sec_simul}

\subsection{Practical implementation of the adaptive choice of $\tau$}
\label{secsecsimul} From the proof of Theorem \ref{maintheorem}, we have the representation
$$\hat{\theta}-\theta_0=-\frac{1}{n}\sum_{i=1}^n V_{\tau}^{-1}\psi(Z_i,\delta_i,X_i;f_1\mathds{1}_{A_{\tau}})+o_P(n^{-1/2}).$$
As in \cite{Stute95}, the function $\psi$ of Theorem \ref{km} can be estimated from the data in the following way by
$$\hat{\psi}(Z,\delta,X;\hat{f_1}\mathds{1}_{A_{\tau}})=\frac{\delta\hat{f_1}(X,Z)}{1-\hat{G}(Z-)}+\int \frac{\int_{y}^{\tau_0}\int_{\mathcal{X}} \hat{f_1}(x,t)d\hat{F}(x,t) dM^{\hat{G}}(y)}{1-\hat{H}(y)},$$
where $\hat{f}_1$ is our kernel estimator of $f_1$ and $\hat{H}$ is the empirical estimator of $H$.
To consistently estimate $\Delta(f_1),$
we use the general technique proposed by \cite{Stute96}, that is
\begin{equation}
\label{varianceestimee}
\hat{\Delta}_{\tau}(f_1) = \frac{1}{n}\sum_{i=1}^n \left[\hat{\psi}(Z_i,\delta_i,X_i;\hat{f}_1)-\frac{1}{n}\sum_{i=1}^n\hat{\psi}(Z_i,\delta_i,X_i;\hat{f}_1)\right]^{\otimes 2},
\end{equation}
where $\otimes 2$ denotes the product of the matrix with its transpose.
A consistent estimator of $V_{\tau}$ can then be computed as
$$\hat{V}_{\tau}=\int \hat{f}^{{ h,\tau}^{-2}}_{\hat{\theta}}(y,\hat{\theta}'x)\hat{J}_0(x,c)\nabla_{\theta}\hat{f}_{\hat{\theta}}^{{ h,\tau}}(y,x)\nabla_{\theta}\hat{f}_{\hat{\theta}}^{{ h,\tau}}(y,x)'
\mathds{1}_{y\in A_\tau}d\hat{F}(x,y).$$
To estimate the asymptotic mean squared error we use
$$\hat{E}^2_{\tau}=\frac{1}{n}\hat{W}_{n,\tau}'\hat{V_{\tau}}^{-1}\hat{V_{\tau}}^{-1}\hat{W}_{n,\tau}.$$

\subsection{Simulation study}

In order to check the finite sample behavior of our estimators of $\theta_0$, we conducted some simulations using a similar model as the one in Delecroix {\it et al.} (2003).
We considered the following regression model,
\[Y_i = \theta_0'X_i+\varepsilon_i, \quad i=1, \ldots,n\]
where $Y_i\in \mathbb{R}$, $\theta_0=(1,0.5,1.4,0.2)'$ and $X_i \sim \otimes^4  \{0.2\,\mathcal{N}(0,1)+0.8\,\mathcal{N}(0.25,2)\}.$
The errors are centered and normally distributed with conditional variance equal to $|\theta_0'X|$.
We used the kernel
 \[K(u) = 2k(u)-k*k(u)\]
where $*$ denotes the convolution product and
\[k(u)=\frac{3}{4}(1-u^2)\mathds{1}_{|u|\leq 1}\]
is the classical Epanechnikov kernel. The censoring distribution was selected to be exponential with parameter $\lambda$ which allows us to fix the proportion of censored responses ($p=25\%$ and $p=40\%$ in our simulations). $\hat{h}$ was chosen using a regular grid between 1 and 1.5.

Our estimator $\hat{\theta}^{\hat \tau}$ was compared with two other estimators, that is $\hat{\theta}^{\infty}$ which does not rely on an adaptive choice of $\tau,$ and $\hat{\theta}^{ADE}$ which is obtained using the average derivative method of \cite{Lu05}.
In the tables below we report our results over $100$ simulations from samples of size $100$ and $200$ for two different rates of censoring. Recalling that the first component of $\theta_0$ is imposed to be one, we only have to estimate the three other components. For each estimator, the Mean Squared Error $ E(\|\hat{\theta}-\theta_0\|^2)$ is decomposed into bias and covariance.

\begin{table}[h]
\begin{center}
\begin{tabular}{|l|c|c|c|}
\hline $p=25\% , n=100$ & Bias & Variance & MSE\\
\hline \hline $\hat{\theta}^{ADE}$ & $\left(\begin{array}{c}
-0.112\\
-0.551\\
-0.155 \end{array}\right)$
 & $\left(\begin{array}{ccc}
0.14 & 0.005 & -0.022 \\
0.005 & 0.075 & 0.016 \\
-0.022 & 0.016 & 0.116
\end{array}\right)$ & 0.6714181 \\
\hline  $\hat{\theta}^{\infty}$ & $\left(\begin{array}{c}
0.057\\
0.215\\
0.048 \end{array}\right)$ & $\left(\begin{array}{ccc}
0.033 & 0.012 & 0.001 \\
0.012 & 0.073 & -0.004 \\
0.001 & -0.004 & 0.027
 \end{array}\right)$ & 0.1841227 \\
\hline $\hat{\theta}^{\hat \tau}$ & $\left(\begin{array}{c}
0.07\\
0.221\\
0.028\\
\end{array}\right)$ & $\left(\begin{array}{ccc}
0.034 & 0.002 & 0.002 \\
0.002 & 0.074 & 0 \\
0.002 & 0 & 0.02
 \end{array}\right)$ & 0.1825980 \\
\hline
\end{tabular}
\end{center}
\caption{Biases, variances and mean squared errors for $25 \%$ of censoring and sampling of size $100$.}
\end{table}

\begin{table}[H]
\begin{center}
\begin{tabular}{|l|c|c|c|}
\hline $p=40\% , n=100$ & Bias & Variance & MSE \\
\hline \hline $\hat{\theta}^{ADE}$ & $\left(\begin{array}{c}
-0.334\\
-0.743\\
-0.158 \end{array}\right)$
 & $\left(\begin{array}{ccc}
0.159 & 0.009 & -0.014 \\
0.009 & 0.268 & 0.048 \\
-0.014 & 0.048 & 0.165
\end{array}\right)$ & 1.280163 \\
\hline  $\hat{\theta}^{\infty}$ & $\left(\begin{array}{c}
0.127\\
0.296\\
0.096 \end{array}\right)$ & $\left(\begin{array}{ccc}
0.11 & -0.034 & -0.01 \\
-0.034 & 0.101 & 0.021 \\
-0.01 & 0.021 & 0.059
 \end{array}\right)$ & 0.3829797 \\
\hline $\hat{\theta}^{\hat \tau}$ & $\left(\begin{array}{c}
0.074\\
0.176\\
0.061\\
\end{array}\right)$ & $\left(\begin{array}{ccc}
0.064 & -0.005 & -0.004 \\
-0.005 & 0.051 & 0.014 \\
-0.004 & 0.014 & 0.069
 \end{array}\right)$ & 0.2239023 \\
\hline
\end{tabular}
\end{center}
\caption{Biases, variances and mean squared errors for $40 \%$ of censoring and sampling of size $100$.}
\end{table}

\begin{table}[H]
\begin{center}
\begin{tabular}{|l|c|c|c|}
\hline $p=25\%,n=200$ &  Bias & Variance & MSE \\
\hline \hline $\hat{\theta}^{ADE}$ & $\left(\begin{array}{c}
-0.189\\
-0.578\\
-0.133 \end{array}\right)$
 & $\left(\begin{array}{ccc}
0.096 & 0.003 & 0.006 \\
0.003 & 0.148 & -0.016 \\
0.006 & -0.016 & 0.131
\end{array}\right)$ & 0.7620268 \\
\hline  $\hat{\theta}^{\infty}$ & $\left(\begin{array}{c}
0.073\\
0.133\\
0.015 \end{array}\right)$ & $\left(\begin{array}{ccc}
0.033 & 0.004 & -0.004 \\
0.004 & 0.023 & 0.002 \\
-0.004 & 0.002 & 0.012
 \end{array}\right)$ & 0.0910719 \\
\hline $\hat{\theta}^{\hat \tau}$ & $\left(\begin{array}{c}
0.034\\
0.107\\
0.014\\
\end{array}\right)$ & $\left(\begin{array}{ccc}
0.007 & 0.001 & 0.004 \\
0.001 & 0.011 & 0 \\
0 & 0 & 0.006
 \end{array}\right)$ & 0.0364064 \\
\hline
\end{tabular}
\end{center}
\caption{Biases, variances and mean squared errors for $25 \%$ of censoring and sampling of size $200$.}
\end{table}

\begin{table}[H]
\begin{center}
\begin{tabular}{|l|c|c|c|}
\hline $p=40\%,n=200$ & Bias & Variance & MSE \\
\hline \hline $\hat{\theta}^{ADE}$ & $\left(\begin{array}{c}
-0.109\\
-0.763\\
-0.053 \end{array}\right)$
 & $\left(\begin{array}{ccc}
0.146 & -0.02 & 0.056 \\
-0.02 & 0.143 & -0.014 \\
0.056 & -0.014 & 0.192
\end{array}\right)$ & 1.078027 \\
\hline  $\hat{\theta}^{\infty}$ & $\left(\begin{array}{c}
0.104\\
0.151\\
0.077 \end{array}\right)$ & $\left(\begin{array}{ccc}
0.109 & 0.008 & 0.042 \\
0.008 & 0.049 & 0.003 \\
0.042 & 0.003 & 0.055
 \end{array}\right)$ & 0.2521227 \\
\hline $\hat{\theta}^{\hat \tau}$ & $\left(\begin{array}{c}
0.043\\
0.14\\
0.021\\
\end{array}\right)$ & $\left(\begin{array}{ccc}
0.018 & -0.001 & 0.002 \\
-0.001 & 0.022 & 0.002 \\
0.002 & 0.002 & 0.014
 \end{array}\right)$ & 0.07533921 \\
\hline  
\end{tabular}
\end{center}
\caption{Biases, variances and mean squared errors for $40 \%$ of censoring and sampling of size $200$.}
\end{table}

To give a precise idea of the number of observations which are removed from the study by choosing $\tau$ adaptively, introduce
$N=\sharp\{1\leq i \leq n, Z_i\leq \hat{\tau}\}.$ In the following table \ref{poi}, we evaluated $E[N]$ in the different cases we considered in the simulation study. We also mention the average weight allocated to the largest (uncensored) data point, first in the case where we consider the whole data set (we denote it Weight$^{\infty}$), then in the truncated data set where we removed all data points with $Z_i\geq \hat{\tau}$  (we denote it Weight$^{\hat \tau}$).

\begin{table}[ht]
\begin{center}
\begin{tabular}{|l|c|c|c|}
\hline& $\hat{\mathbb E}(N)$ & Weight$^{\infty}$ & Weight$^{\hat \tau}$ \\
\hline \hline $n=100,p=25\%$ & $90$
 & $0.0667$ & $0.0204$ \\
\hline  $n=100,p=40\%$ & 87 & $0.124$ & $0.0236$ \\
\hline $n=200,p=25\%$ & $185$ & $0.0402$ & $0.0119$  \\
\hline  $n=200,p=40\%$ & 172 & $0.0997$ & $0.0122$ \\
\hline 
\end{tabular}
\label{poids}
\end{center}
\caption{Last observed data in the truncating model and weight allocated to the largest observation in each model for different sample sizes and censoring rates.}
\label{poi}
\end{table}


Clearly the MSE deteriorates when the percentage of censoring increases. According to the simulations, $\hat{\theta}^{\hat \tau}$ and $\hat{\theta}^{\infty}$ outperform $\hat{\theta}^{ADE},$ while, as expected, choosing adaptively $\tau$ improves the quality of the estimation. This is not obvious in the case where there are only $25\%$ of censoring. However, in the case where the level of censoring is high, estimation of the tail of the distribution by Kaplan-Meier estimators becomes more erratic, and the importance of choosing a proper truncation appears in the significant difference between the MSE of $\hat{\theta}^{\hat \tau}$ and $\hat{\theta}^{\infty}.$ Moreover, the importance of truncation becomes obvious if we look at table \ref{poi}. We see that, in the case where there is 40 \% of censoring, the weight allocated to the largest data-point if we do not use truncation can be up to 10 times (approximatively) the weight allocated to the largest observation in the truncated data set. The ratio is less important in the case where there is 25 \% of censoring, but still consequent (in this case, the ratio is approximatively 3). Therefore, it seems that,  considering the whole data set, the weight allocated to the largest observation can have a too  strong  influence on the estimation procedure, which explains the difference of performance of the estimators with or without truncation.

\subsection{Example : Stanford Heart Transplant Data}

We now illustrate our method using data from the Stanford Heart
Transplant program. This data set was initially studied by \cite{Miller82}. 184 of 249 patients in this program
received a heart transplantation between October 1967 and February
1980. From this data, we considered the survival time as the
response variable $Z$, age as the first component of $X$ and the
square of age as the second {component}. Patients alive beyond February 1980 were considered censored. For easier comparison to previous
work on this data set, we concentrate our analysis on the 157
patients out of 184 who had complete tissue typing. Among these 157 cases, 55 were censored.

Several methods of estimation have already been applied to this
data set to estimate the following regression model,
\begin{equation}
Z=\alpha+\beta'X+\varepsilon(X), \label{modelheart} \end{equation}
where $\beta=(\beta_1,\beta_2)',$ $E[\varepsilon(X)|X]=0,$ see \cite{Miller82}, Wei {\it et al.} (1990), Stute {\it et al.} (2000). Furthermore, nonparametric
lack-of-fit tests have shown that the regression model
(\ref{modelheart}) seemed reasonable, see Stute {\it et al.} (2000) and \cite{Lopez08}.
Therefore it seems to us appropriate to experiment our model on
this data set. This strengthens the assumption on the residual, by
assuming that $\varepsilon(X)=\varepsilon(\theta_0'X),$ where
$\theta_0=(1,\beta_2/\beta_1)',$ but allows more flexibility on the
regression function.

In the following table, we present our estimators and recall the
values of the estimators of $\beta_2/\beta_1$ for the linear
regression model (\ref{modelheart}). We first computed
$\hat{\theta}^{\infty},$ which is our estimator using the whole data set,
that is with $\tau=+\infty,$ and compared it to the one obtained
by choosing $\tau$ from the data as in section \ref{secsecsimul}.
In this last case, $\hat{\tau}=Z_{(90)}$ where $Z_{(i)}$ denotes
the $i-$th order statistic, this means that it conducted us to
remove the $67$ largest observations to estimate $\theta_0$ (but
not to estimate Kaplan-Meier weights, which were computed using
the whole data set). We computed Weight$^{\infty}=0.0397,$ and Weight$^{\hat{\tau}}= 0.0076$ for the truncated data. Adaptive bandwidth was $1.7$ for
$\hat{\theta}^{\infty},$ and $1.3$ for $\hat{\theta}^{\hat \tau}.$ The estimated
value of the mean-squared error was $E^2_{\infty}=0.1089375$
and $E^2_{\hat{\tau}}=0.01212701$ for $\hat{\theta}^{\infty}$ and
$\hat{\theta}^{\hat \tau}$ respectively.

\begin{table}[h]
\begin{center}
\begin{tabular}{|l|c|}
\hline  & Estimator of $\theta_{0,2}=\beta_2/\beta_1$ \\
\hline \hline Miller and Halpern & -0.01588785 \\
\hline Wei \it{et al.} & 63.75 \\
\hline Stute \it{et al.} & -0.01367034 \\
\hline $\hat{\theta}^{\infty}$ (without adaptive choice of $\tau$) & -0.07351351\\
\hline $\hat{\theta}^{\hat \tau}$ (with adaptive choice of $\tau$) & -0.0421508  \\
\hline
\end{tabular}
\end{center}
\caption{Comparison of different estimators of $\theta_{0,2}$.}
\end{table}
Our estimators seem relatively close to the ones obtained by \cite{Miller82} and Stute {\it et al.} (2000) using respectively the Buckley-James method and the Kaplan-Meier integrals method for the linear regression model.
%

\section{Proof of Main Lemma}\label{lemme}
First, the same arguments as in Delecroix {\it et al.} (2006) apply to replace $\hat{J}_0$ by $J_0.$ 
Define $J_{\theta}(x,c)=\mathds{1}_{f_{\theta'X}(\theta'x)\geq c}.$ From
Assumption \ref{lgnu} on the density of $\theta'x$, deduce that,
on shrinking neighbourhoods of $\theta_0,$ $J_0(x,c)$ can be
replaced by $J_{\theta}(x,c/2).$ Using a Taylor expansion, write
\begin{align*}
L^{\tau}_n(\theta,\hat{f}^{h,\tau},J_0)-L^{\tau}_n(\theta,f^{\tau},J_0) &= \sum_{i=1}^n \delta_i W_{in}\mathds{1}_{Z_i\in A_{\tau}}\log\left(\frac{\hat{f}^{h,\tau}_{\theta}(Z_i,\theta'X_i)}{f_{\theta}^{{\tau}}(Z_i,\theta'X_i)}\right)J_0(X_i,c) \\
&= \sum_{i=1}^n \frac{\delta_i W_{in}\mathds{1}_{Z_i\in A_{\tau}}\left(\hat{f}^{h,\tau}_{\theta}(Z_i,\theta'X_i)-f_{\theta}^{{\tau}}(Z_i,\theta'X_i)\right)J_0(X_i,c)}{f_{\theta}^{{\tau}}(Z_i,\theta'X_i)} \\
& -\sum_{i=1}^n \frac{\delta_i W_{in}\mathds{1}_{Z_i\in A_{\tau}}\left[\hat{f}^{h,\tau}_{\theta}(Z_i,\theta'X_i)-f_{\theta}^{{\tau}}(Z_i,\theta'X_i)\right]^2J_0(X_i,c)}{\phi(f_{\theta}^{{\tau}}(Z_i,\theta'X_i),\hat{f}^{h,\tau}_{\theta}(Z_i,\theta'X_i))^2} \\
&= A^{\tau}_{1n}(\theta,\hat{f}^{h,\tau})-B^{\tau}_{1n}(\theta,\hat{f}^{h,\tau})
\end{align*}
where $\phi(f_{\theta}^{{\tau}}(Z_i,\theta'X_i),\hat{f}^{h,\tau}_{\theta}(Z_i,\theta'X_i))$ is between $\hat{f}^{h,\tau}_{\theta}(Z_i,\theta'X_i)$ and $f_{\theta}^{{\tau}}(Z_i,\theta'X_i).$

\textbf{Step 1.} We first study $A_{1n}.$ A Taylor expansion leads to the following decomposition,
\begin{align*}
A^{\tau}_{1n}&= (\theta-\theta_0)'\sum_{i=1}^n \frac{\delta_i W_{in}\mathds{1}_{Z_i\in A_{\tau}}\big(\nabla_{\theta}\hat{f}^{h,\tau}_{\theta_0}(Z_i,X_i)-\nabla_{\theta}f_{\theta_0}^{{\tau}}(Z_i,X_i)\big)J_{\theta}(X_i,c/2)}{f_{\theta}^{{\tau}}(Z_i,\theta'X_i)} \\ & \quad+ (\theta-\theta_0)'\left[\sum_{i=1}^n \frac{\delta_i W_{in}\mathds{1}_{Z_i\in A_{\tau}}\big(\nabla^2_{\theta}\hat{f}^{h,\tau}_{\tilde{\theta}}(Z_i,X_i)-\nabla^2_{\theta}f_{\tilde{\theta}}^{{\tau}}(Z_i,X_i)\big)J_{\theta}(X_i,c/2)}{{2}f_{\theta}^{{\tau}}(Z_i,\theta'X_i)}\right](\theta-\theta_0) \\
& \quad+ \frac{1}{n}\sum_{i=1}^n \frac{\delta_iW_{in}\mathds{1}_{Z_i\in A_{\tau}}\big(\hat{f}^{h,\tau}_{\theta_0}(Z_i,\theta_0'X_i)-f_{\theta_0}^{{\tau}}(Z_i,\theta_0'X_i)\big)}{f_{\theta}^{{\tau}}(Z_i,\theta'X_i)f_{\theta_0}^{{\tau}}(Z_i,\theta_0'X_i)}\\
& \quad \times (f_{\theta_0}^{{\tau}}(Z_i,\theta_0'X_i)-f_{\theta}^{{\tau}}(Z_i,\theta'X_i))J_0(X_i,c)J_{\theta}(X_i,c/2)+
{A^{\tau}_{1n}(\theta_0,\hat{f}^{h,\tau})}\\
&= {A^{\tau}_{1n}(\theta_0,\hat{f}^{h,\tau})}+(\theta-\theta_0)'A^{\tau}_{2n}(\theta_0,\hat{f}^{h,\tau})+(\theta-\theta_0)'A^{\tau}_{3n}(\tilde{\theta},\hat{f}^{h,\tau})(\theta-\theta_0)+A^{\tau}_{4n}(\theta,\hat{f}^{h,\tau}),
\end{align*}
for some $\tilde{\theta}$ between $\theta$ and $\theta_0.$ Observe that, using the uniform consistency of $\nabla_{\theta}^2\hat{f}^{h,\tau}_{\theta}$ (deduced from Proposition \ref{convergence_rate} and Lemma \ref{leplusdur}), we obtain $\sup_{\tilde{\theta}\in \Theta_n,\tau\leq \tau_0,h\in \mathcal{H}_n}A^{\tau}_{3n}(\tilde{\theta},\hat{f}^{h,\tau})=o_P(1).$ We now study $A^{\tau}_{2n}(\theta_0,\hat{f}^{h,\tau}).$ Using the expression (\ref{explicit_jump}) of the jumps of Kaplan-Meier estimator, observe that
\begin{align*}
& A^{\tau}_{2n}(\theta,\hat{f}^{h,\tau})
\\ &\quad= \sum_{i=1}^n \frac{W_i^*\mathds{1}_{Z_i\in A_{\tau}}\big(\nabla_{\theta}\hat{f}^{h,\tau}_{\theta_0}(Z_i,X_i)-\nabla_{\theta}f_{\theta_0}^{{\tau}}(Z_i,X_i)\big)J_{\theta}(X_i,c/2)}{f_{\theta}^{{\tau}}(Z_i,\theta'X_i)}\\ &\qquad + \frac{1}{n}\sum_{i=1}^nW_i^*Z_G(Z_i-)\frac{\delta_i\mathds{1}_{Z_i\in A_{\tau}}\big(\nabla_{\theta}\hat{f}^{h,\tau}_{\theta_0}(Z_i,X_i)-\nabla_{\theta}f_{\theta_0}^{{\tau}}(Z_i,X_i)\big)J_{\theta}(X_i,c/2)}{f_{\theta}^{{\tau}}(Z_i,\theta'X_i)}
\\
&\quad= A^{\tau}_{21n}(\theta,\hat{f}^{h,\tau})+A^{\tau}_{22n}(\theta,\hat{f}^{h,\tau}),
\end{align*}
where \[Z_G(t)=\frac{\hat{G}(t)-G(t)}{1-\hat{G}(t)}.\]
The term $A^{\tau}_{22n}$ can be bounded using (\ref{gill}), (\ref{zhou}) and Lemma \ref{leplusdur}, by
$$\sup_{\tau \leq \tau_0}|A_{22n}(\theta,\hat{f}^{h,\tau})|\leq o_P(n^{-1/2})\times n^{-1}\sum_{i=1}^n \delta_i[1-G(Z_i-)]^{-1},$$
and the last term is $O_P(1)$ since it has finite expectation. Now for $A^{\tau}_{21n},$ first replace $\theta$ at the denominator by $\theta_0.$ We have
\begin{align*}
A^{\tau}_{21n}(\theta,\hat{f}^{h,\tau}) & =\sum_{i=1}^n \frac{W_i^*\mathds{1}_{Z_i\in A_{\tau}}(\nabla_{\theta}\hat{f}^{h,\tau}_{\theta_0}(Z_i,X_i)-\nabla_{\theta}f_{\theta_0}^{{\tau}}(Z_i,X_i))J_{0}(X_i,c/4)}
{f_{\theta_0}^{{\tau}}(Z_i,\theta_0'X_i)}\\
& \quad +R^{\tau}_{n}(\theta,h)(\theta-\theta_0),
\end{align*}
with $\sup_{\theta\in \Theta_n,\tau\leq \tau_0,h\in \mathcal{H}_n}|R_n^{\tau}(\theta,h)|=o_P(1)$ from Assumption \ref{lgnu} and the uniform consistency of $\nabla_{\theta}\hat{f}^{h,\tau}_{\theta_0}$ deduced from Proposition \ref{convergence_rate} and Lemma \ref{leplusdur}.
Then use Assumption \ref{donsker1} and Proposition \ref{classe de Donsker}. Using the equicontinuity property of Donsker classes (see e.g. \cite{Vaart96} or \cite{Vaart98}), we obtain that
\begin{align*}
A^{\tau}_{2n}(\theta,\hat{f}^{h,\tau}) & =\iint \frac{\left[\nabla_{\theta}\hat{f}^{h,\tau}_{\theta_0}(y,x)-\nabla_{\theta}f_{\theta_0}^{{\tau}}(y,x)\right]\mathds{1}_{y\in A_{\tau}}J_{0}({ x},c/4)d\mathbb{P}(x,y)}{f_{\theta_0}^{{\tau}}(y,u)}\\
& \quad +o_P(n^{-1/2}),
\end{align*}
where the $o_P-$rate does not depend on $\theta,$ $h,$ nor $\tau.$
From classical kernel arguments, \linebreak $\sup_{y,x,\tau}|\int\big(\nabla_{\theta}f^{*h,\tau}_{\theta_0}(y,x)-\nabla_{\theta}
f^{{\tau}}_{\theta_0}(y,x)\big)\mathds{1}_{y\in
A_{\tau}}J_{0}({ x},c/4)d\mathbb{P}(x,y)|=O_{\mathbb{P}}(h^4)=o_\mathbb{P}(n^{-1/2}),$ since $nh^8\rightarrow 0.$ Then,
Lemma \ref{leplusdur2} concludes the proof for
$A_{2n}^{\tau}(\theta,\hat{f}^{h,\tau})$. $A^{\tau}_{4n}(\theta,\hat{f}^{h,\tau})$ can be handled similarly.
%

%

\textbf{Step 2.} $B^{\tau}_{1n}$ can be rewritten as
\begin{align*}
&B^{\tau}_{1n}(\theta,\hat{f}^{h,\tau})\\
&\quad=\sum_{i=1}^n \delta_iW_{in}\mathds{1}_{Z_i\in A_{\tau}}J_{\theta}(X_i,c/2)\frac{\{(\theta-\theta_0)'[\nabla_{\theta}\hat{f}^{h,\tau}_{\theta_0}(Z_i,X_i)-\nabla_{\theta}f_{\theta_0}^{{\tau}}(Z_i,X_i)]\}^2}{\phi(f_{\theta}^{{\tau}}(Z_i,\theta'X_i),\hat{f}_{\theta}^{{h,\tau}}(Z_i,\theta'X_i))^2}
\\
&\qquad+ 2\sum_{i=1}^n \delta_iW_{in}J_{\theta}(X_i,c/2)\mathds{1}_{Z_i\in A_{\tau}}[\hat{f}^{h,\tau}_{\theta_0}(Z_i,\theta_0'X_i)-f_{\theta_0}^{{\tau}}(Z_i,\theta_0'X_i)]\\
&\qquad \times (\theta-\theta_0)'[\nabla_{\theta}\hat{f}^{h,\tau}_{\tilde{\theta}}(Z_i,X_i)-\nabla_{\theta}f_{\tilde{\theta}}^{{\tau}}(Z_i,X_i)][\phi(f_{\theta}^{{\tau}}(Z_i,\theta'X_i),\hat{f}_{\theta}^{{h,\tau}}(Z_i,\theta'X_i))^2]^{-1}
\\
&\qquad+B^{\tau}_{4n}(\theta_0,\hat{f}^{h,\tau})+o_P(\|\theta-\theta_0\|^2)& \\ &\quad=(\theta-\theta_0)'B^{\tau}_{2n}(\theta_0,\hat{f}^{h,\tau})(\theta-\theta_0)+(\theta-\theta_0)'B^{\tau}_{3n}(\theta,\hat{f}^{h,\tau})+B^{\tau}_{4n}(\theta_0,\hat{f}^{h,\tau})+o_P(\|\theta-\theta_0\|^2)
\end{align*}
for some $\tilde{\theta}$ between $\theta$ and $\theta_0$. The third term does not depend on $\theta.$ For $B_{2n}^{\tau},$ use the uniform consistency of $\nabla_{\theta_0}\hat{f}^{h,\tau}_{\theta_0}$ (Proposition \ref{convergence_rate} and Lemma \ref{leplusdur}) to obtain $\sup_{\tau\leq \tau_0,h\in \mathcal{H}_n}|B^{\tau}_{2n}(\theta,\hat{f}^{h,\tau})|=o_P(n^{-1/2}).$ Finally, for $B^{\tau}_{3n}(\theta,\hat{f}^{h,\tau}),$ from a Taylor expansion,
\begin{align*}
B^{\tau}_{3n}(\theta,\hat{f}^{h,\tau}) & =
2\sum_{i=1}^n \frac{\delta_iW_{in}\mathds{1}_{Z_i\in A_{\tau}}J_{\theta}(X_i,c/2)}{\phi(f_{\theta}^{{\tau}}(Z_i,\theta'X_i),\hat{f}^{h,\tau}_{\theta}(Z_i,\theta'X_i))^2}\\
& \quad \times[\hat{f}^{h,\tau}_{\theta_0}(Z_i,\theta_0'X_i)-f_{\theta_0}^{{\tau}}(Z_i,\theta_0'X_i)][\nabla_{\theta}\hat{f}^{h,\tau}_{\theta_0}(Z_i,X_i)-\nabla_{\theta}f_{\theta_0}^{{\tau}}(Z_i,X_i)] \\
&\quad +(\theta-\theta_0)'R_n^{\tau}(\theta,\hat{f}^{h,\tau}),\end{align*}
with $\sup_{\theta\in \Theta_n,\tau\leq \tau_0,h\in
\mathcal{H}_n}R_n^{\tau}(\theta,\hat{f}^{h,\tau})=o_P(1),$ from
Proposition \ref{convergence_rate} and Lemma \ref{leplusdur}. For
the main term, the product of the uniform convergence rates of
$\hat{f}^{h,\tau}_{\theta_0}$ and
$\nabla_{\theta}\hat{f}^{h,\tau}_{\theta_0}$ obtained from
Proposition \ref{convergence_rate} and {Lemma \ref{leplusdur}} is
$o_P(n^{-1/2})$ for $h\in \mathcal{H}_n$.

\section{Conclusion}

We proposed a new estimation procedure of a conditional density under a single-index assumption and random censoring. This procedure is an extension of the approach of Delecroix {\it et al.} (2003) in the case of a censored response. One of the advantage of this model is that it relies on fewer assumptions as a Cox regression model, in the case where the random variables of the model are absolutely continuous. By showing that estimating in this semiparametric model is asymptotically equivalent to estimating in a parametric one (unknown in practice), we obtain a $n^{-1/2}-$rate for the estimator of the index. This estimator can then be used to estimate the conditional density or the conditional distribution function by using traditional nonparametric estimator under censoring. A new feature of our procedure, is that it provides an adaptively driven choice of the bandwidth involved in the kernel estimators we used, and that it also provides an adaptive choice of a truncation parameter needed to avoid problems caused by the bad behavior of Kaplan-Meier {estimators} in the tail of the distribution. In this specific problem, this truncation does not introduce some additional bias in the procedure, and seems, according to our simulations, to increase the quality of the estimator, especially in the case where the proportion of censored responses is important. Our way of choosing $\tau$ was motivated by minimizing the MSE in the estimation of $\hat{\theta}.$ However, our method could be easily adapted to other kinds of criteria which, for example more focus on the error in estimating one specific direction, or on the error in the estimation of the conditional density itself.

\section{Appendix}\label{technical}

\subsection{Kaplan-Meier integrals for the parametric case}
We first recall a classical asymptotic representation of integrals with respect to $\hat{F}.$ See \cite{Stute95}, \cite{Stute96} and S{\'a}nchez Sellero {\it et al.} (2005).

\begin{theorem}
\label{km}
Let $\mathcal{F}$ be a $VC-$class of functions with envelope $\Phi$ such as
\begin{equation}
\label{nullity}
\Phi(x,y)= 0, \; for \;  all \; y\geq \tau_0,
\end{equation}
where $\tau_0\leq \tau_H.$ We have the following asymptotic i.i.d. representation, for all $\phi\in \mathcal{F},$
$$ \int \phi(x,y)d\hat{F}(x,y) = \frac{1}{n}\sum_{i=1}^n \psi(Z_i,\delta_i,X_i;\phi)+R(\phi),$$
where $\sup_{\phi \in \mathcal{F}}|R(\phi)|=O_{a.s.}([\log n]^3n^{-1}),$ and
$$\psi(Z_i,\delta_i,X_i;\phi)=\frac{\delta\phi(X_i,Z_i)}{1-G(Z_i-)}+\int \frac{\int_{y}^{\tau_0}\int_{\mathcal{X}} \phi(x,t)dF(x,t) dM_i^G(y)}{1-H(y)},$$
where $M_i^G(y)=(1-\delta_i)\mathds{1}_{Z_i\leq y}-\int_{-\infty}^y \mathds{1}_{Z_i\geq t}[1-G(t-)]^{-1}dG(t)$ is a martingale with respect to the filtration $\mathcal{G}_y=\{(Z_i,\delta_i,X_i)\mathds{1}_{Z_i\leq y}\}.$
Define $\Delta(\phi)=Var(\psi(Z,\delta,X;\phi)).$ Then it follows that
$$\sqrt{n}\int \phi(x,y)d[\hat{F}-F](x,y)\Longrightarrow \mathcal{N}(0,\Delta(\phi)).$$
\end{theorem}

Initially, the result of Stute was derived for a single function $\phi.$ Furthermore, Theorem 1.1 in \cite{Stute96} gives a convergence rate which is only $o_P(n^{-1/2})$ for the remainder term, however an higher convergence rate is obtained in his proof of Theorem 1.1 for functions satisfying (\ref{nullity}), which is the only case considered in our work. To obtain uniformity on a $VC-$class of functions, see S{\'a}nchez Sellero {\it et al.} (2005) who provided a more general representation that extends the one of Stute in the case where $Y$ is right-censored and left-truncated. Their result is really useful since it provides, as a corollary, uniform law of large numbers results and uniform central limit theorem. The representation we present in our Theorem \ref{km} is a simple rewriting of Stute's representation. Theorem \ref{km} is then a key ingredient to prove Lemma \ref{parametric}.

\begin{proof}[Proof of Lemma \ref{parametric}]
We directly show the second part of the Lemma, since the first can be studied from similar techniques. From a Taylor expansion,
\begin{align}
L^{\tau}_n(\theta,f^{{\tau}},J_0)&= (\theta-\theta_0)'\sum_{i=1}^n \delta_iW_{in}J_0(X_i,c)\mathds{1}_{Z_i\in A_{\tau}}\frac{\nabla_{\theta}f_{\theta_0}^{{\tau}}(Z_i,X_i)}{f_{\theta_0}^{{\tau}}(Z_i,\theta_0'X_i)}  \nonumber \\
&\quad +\frac{1}{2}(\theta-\theta_0)'\sum_{i=1}^n \delta_iW_{in}J_0(X_i,c)\mathds{1}_{Z_i\in A_{\tau}}\nabla^2_{\theta}[\log f_{\tilde{\theta}}^{{\tau}}](Z_i,X_i)(\theta-\theta_0)\nonumber \\
&\quad +T_{4n}(\theta_0), \label{development}
\end{align}
for some $\tilde{\theta}$ between $\theta_0$ and $\theta.$
Theorem \ref{km} provides an i.i.d. representation for the first term (which corresponds to $W_{n,\tau}$ in Lemma \ref{parametric}), while, from Assumption \ref{regularity}, the family of functions $\nabla^2_{\theta}[\log f_{\tilde{\theta}}^{{\tau}}](y,x)\mathds{1}_{y\in A_{\tau}}$ is a $VC-$class of functions satisfying (\ref{nullity}). Hence the sum in the second term of (\ref{development}) tends to $V$ almost surely using an uniform law of large numbers property.
\end{proof}

\subsection{The gradient of $f$}
In the following for any function $\varphi$ we will denote by $\varphi_h^{(n)}(\cdot)$ the expression $h^{-n}\varphi^{(n)}(\cdot /h)$ such as, for example $K'_h(\cdot) =h^{-1}K'\left(\cdot/h\right).$

\begin{proposition}
\label{expectation} Let
\[f'_{\tau}(y,u)=\partial_u f^{\tau}_{\theta_0}(y,u).\]
We have
$$\nabla_{\theta}f^{\tau}_{\theta_0}(y',x)=xf_{1,y,\tau}(y',\theta_0'x)+f_{2,y,\tau}(y',\theta_0'x),$$
with
\begin{align*}
f_{1,\tau}(y,\theta_0'x) &= f'_{\tau}(y,\theta_0'x), \\
f_{2,\tau}(y,\theta_0'x) &= -f'_{\tau}(y,\theta_0'x)E\left[X|\theta_0'X\right].
\end{align*}
In particular, $E[\nabla_{\theta}f_{\theta_0}^{\tau}(Y,X)|\theta_0'X]=0.$
\end{proposition}

\begin{proof}
Direct adaptation of Lemma 5A in \cite{Dominitz05}.
\end{proof}

\subsection{Convergence properties of $f^{*h,\tau}$}

We first recall some classical properties on kernel estimators. Consider the class of functions $\mathcal K$ introduced in Assumption \ref{noyau}. Let $N(\varepsilon, \mathcal K, d_Q)$ be the minimal number of balls $\{g:d_Q(g,g')<\varepsilon\}$ of $d_Q$-radius $\varepsilon$ needed to cover $\mathcal K$. For $\varepsilon>0$, let $N(\varepsilon,\mathcal K)=\sup_QN(\kappa\varepsilon,\mathcal{K},d_Q)$, where the supremum is taken over all probability measures $Q$ on $(\mathbb{R}^d,\mathcal B)$, $d_Q$ is the $L_2(Q)$-metric. From \cite{Nolan87}, it can easily be seen that, using a kernel $K$ satis\-fying Assumption \ref{noyau}, for some $C>0$ and ${\nu>0, N(\varepsilon, \mathcal K)\leq C\varepsilon^{-\nu}, 0<\varepsilon<1}$.

\begin{proposition}
\label{convergence_rate} Under assumption \ref{noyau} we have, for some $c>0$
\begin{align}\label{uniff}
\sup_{x,y,h,\tau}\left|f^{*h,\tau}_{\theta_0}(y,\theta_0'x)-f^{\tau}_{\theta_0}(y,\theta_0'x)\right|\mathds{1}_{y\in A_\tau}J_0(x,c)
& = O_P\left(n^{-1/2}h^{-1}[\log n]^{1/2}\right),\\
\label{unifgrad}
\sup_{x,y,h,\tau}\left|\nabla_{\theta}f^{*h,\tau}_{\theta_0}(y,x)-\nabla_{\theta}f^{\tau}_{\theta_0}(y,x)\right|\mathds{1}_{y\in
A_\tau}J_0(x,c)
& = O_P\left(n^{-1/2}h^{-2}[\log n]^{1/2}\right),\\
\sup_{x,y,h,\tau,\theta}\left|\nabla^2_{\theta}f^{*h,\tau}_{\theta}(y,x)-\nabla^2_{\theta}f^{\tau}_{\theta}(y,x)\right|\mathds{1}_{y\in A_\tau}J_{\theta}(x,c)
& = O_P\left(n^{-1/2}h^{-3}[\log n]^{1/2}\right).\label{unif2}
\end{align}
\end{proposition}

\begin{proof}
(\ref{uniff}) is a direct application of Theorems 1 and 4 in \cite{Einmahl05}. For (\ref{unifgrad}), we only show the convergence for the term
$$\hat{r}^{h,\tau}_{\theta_0}(x,y):=\frac{1}{h}\sum_{i=1}^n\delta_iW_i^*\mathds{1}_{Z_i\in A_{\tau}}J_0(x,c)(X_i-x)K'_h(X_i'\theta_0-x'\theta_0)K_h(Z_i-y).$$
Define
$$\bar{r}^{h,\tau}_{\theta_0}(x,y) =\frac{1}{h}\mathbb{E}\,\big[\mathds{1}_{Y\in A_{\tau}}J_0(x,c)(X-x)K'_h(X'\theta_0-x'\theta_0)K_h(Y-y)\big]$$
and
$$r^{\tau}_{\theta_0}(x,y)=\left.\frac{\partial}{\partial u}\left\{\mathbb{E}\,\big[(X-x)|\theta_0'X=u,Y=y\big]\mathds{1}_{y\in A_{\tau}}J_0(x,c)f_{\theta_0'X,Y}(u,y)\right\}\right|_{u=\theta_0'x}.$$
Note that, from our assumptions $r^{\tau}_{\theta_0}$ is a finite quantity.
Next, Theorem 4 in \cite{Einmahl05} yields :
$$\sup_{x,y,h,\tau}\left|\hat{r}^{h,\tau}_{\theta_0}(x,y)-\bar{r}^{h,\tau}_{\theta_0}(x,y)\right|\mathds{1}_{y\in A_{\tau}}= O_P(n^{-1/2}h^{-2}[\log n]^{1/2}).$$
For the bias term, $\sup_{x,y,h,\tau}\left|\bar{r}^{h,\tau}_{\theta_0}(x,y)-r^{\tau}_{\theta_0}(x,y)\right|\mathds{1}_{y\in A_{\tau}}
= O(h^4)=o(n^{-1/2}),$ (see e.g. \cite{Bosq97}). As a consequence,

%

\[\sup_{x,y,h,\tau}\left|\hat{r}^{h,\tau}_{\theta_0}(x,y)-r^{\tau}_{\theta_0}(x,y)\right|\mathds{1}_{y\in A_{\tau}}= O_P(n^{-1/2}h^{-2}[\log n]^{1/2}).\]
For (\ref{unif2}), we also need an uniformity with respect to $\theta.$ The result can be deduced from the uniform convergence (with respect to $\theta,$ $x,$ $u$) of quantities such as
\begin{equation}
\label{formegenerale}
S_n^{h,\tau}(\theta,x,y,\beta)=\frac{1}{h^{2}}\sum_{i=1}^n \delta_iW_i^*\phi(Z_i,X_i,\theta)\nabla^{\beta}_{\theta}K\left(\frac{\theta'X_i-\theta'x}{h}\right)K\left(\frac{Z_i-y}{h}\right),
\end{equation}
where $\nabla^{\beta}_{\theta}K([\theta'X_i-\theta'x]h^{-1})$ for $\beta=1$ (resp. for $\beta=2$) denotes the gradient vector
of function $K([\theta'X_i-\theta'x]/h)$ (resp. Hessian matrix) with respect to $\theta$ and evaluated at $\theta,$ and where $\phi$ is a bounded function with respect to $\theta$ and $x$. The function $\phi$ we consider is $\phi(Z,X,\theta)=f^{\tau}_{\theta'X}(\theta'X)^{-1}\mathds{1}_{{Z\in A_\tau}}J_0(x,c)$ with the convention $0/0=0$ and where $f^{\tau}_{\theta'X}(\theta'X)$ is the conditional density of $\theta'X$ given $Y\in A_\tau.$ (\ref{formegenerale}) can be studied using the same method as \cite{Einmahl05}. For this, observe that the family of functions $\{(X,Z)\rightarrow \nabla_{\theta}^{\beta}K([\theta'X-\theta'x]h^{-1})K([Z-y]h^{-1}),\theta\in \Theta,x,y\}$ satisfies the Assumptions of Proposition 1 in \cite{Einmahl05} (see Lemma 22 (ii) in \cite{Nolan87}). Hence, apply Talagrand's inequality (\cite{Talagrand94}, see also \cite{Einmahl05}) to obtain that
$$\sup_{\theta,x,y,h,\tau}|S_n^{h,\tau}(\theta,x,y,\alpha)-E[S_n^{h,\tau}(\theta,x,y,\alpha)]|\mathds{1}_{y\in A_{\tau}}=O_P(n^{-1/2}[\log n]^{1/2}h^{-1-\beta}).$$ Again, the bias term converges uniformly at rate $O(h^4).$
%
%
%
%
\end{proof}

\subsection{The difference between $f^*$ and $\hat{f}$}

\subsubsection{Convergence rate of $\hat{f}$}

In this section, we show that replacing $f^{*h,\tau}$ by $\hat{f}^{h,\tau}$ (which is the estimator used in practice) does not modify the rate of convergence. To give the intuition of this results, observe that $f^{*h,\tau}$ was obtained from $\hat{f}^{h,\tau}$ by replacing $\hat{G}$ by $G.$ Let us recall some convergence properties of $\hat{G}.$ We have
\begin{align}
\sup_{t\leq \tau_0}|\hat{G}(t)-G(t)| &= O_P(n^{-1/2}), \label{gill} \\
\sup_{t\leq \tau_0}\frac{1-G(t)}{1-\hat{G}(t)} &= O_P(1). \label{zhou}
\end{align}
See \cite{Gill83} for (\ref{gill}) and \cite{Zhou92} for (\ref{zhou}). From (\ref{gill}), we see that the convergence rate of $\hat{G}$ is faster than the convergence rate of $f^{*h,\tau},$ which explains the asymptotic equivalence of $\hat{f}^{h,\tau}$ and $f^{*h,\tau}.$
Lemma \ref{leplusdur}
makes things more precise and also gives a representation of the difference between $\nabla_{\theta}f_{\theta_0}^{*h,\tau}$ and $\nabla_{\theta}\hat{f}_{\theta_0}^{h,\tau}$ which is needed in the proof of Main Lemma. Also required to prove our Main Lemma, Lemma \ref{leplusdur2} below gives a technical result on the integral of this difference.

\begin{lemma}
\label{leplusdur}
Under the Assumption of Lemma \ref{parametric}, we have for some $c>0$
\begin{align}
\label{01} \sup_{x,y,h,\tau}\left|\hat{f}_{\theta_0}^{h,\tau}(y,\theta'x)-f^{*h,\tau}_{\theta_0}(y,\theta'x)\right|\mathds{1}_{y\in A_\tau}J_0(x,c) &= O_P(n^{-1/2}), \\
\label{02} \sup_{x,y,h,\tau}\left|\nabla_{\theta}\hat{f}_{\theta_0}^{h,\tau}(y,x)-\nabla_{\theta}f^{*h,\tau}_{\theta_0}(y,x)\right|\mathds{1}_{y\in A_\tau}J_0(x,c) &= O_P(n^{-1/2}h^{-1}), \\ \label{03}
\sup_{x,y,h,\tau,\theta}\left|\nabla^2_{\theta}\hat{f}_{\theta}^{h,\tau}(y,x)-\nabla^2_{\theta}f^{*h,\tau}_{\theta}(y,x)\right|\mathds{1}_{y\in A_\tau}J_{\theta}(x,c) &= O_P(n^{-1/2}h^{-2}).
\end{align}
Furthermore, for $x$ such as $J_0(x,c)\neq 0$,
\begin{align}
\label{04} \left(\nabla_{\theta}\hat{f}_{\theta_0}^{h,\tau}(y,x)-\nabla_{\theta}f^{*h,\tau}_{\theta_0}(y,x)\right) &= \int \frac{\int_{{\mathcal{X}}}\int_{t}^{\tau_0}g^{h,\tau}_{f,x,y}(x_2,y_2)d\mathbb{P}(x_2,y_2)d\bar{M}^G(t)}{1-H(t)} \nonumber\\
        & \quad +R_n(\tau,h,x,y),
\end{align}
where $\bar{M}^G(y)=n^{-1}\sum_{i=1}^nM_i^G(y),$ $M_i^G$ is defined in Theorem \ref{km}, $\sup_{x,y,\tau,h}|R_n(\tau,h,x,y)|=$ \\ $O_{P}((\log n)^{1/2}n^{-1}h^{-3})$ and $g_{f,x,y}^h$ is defined by
\begin{align*}
g^{h,\tau}_{f,x_1,y_1}(x_2,y_2) &= \frac{1}{h}
\frac{(x_1-x_2)K'_h(\theta_0'x_1-\theta_0'x_2)K_h(y_1-y_2)}{f^{\tau}_{\theta_0'X}(\theta_0'x_1)} \\
& \quad-
\frac{K_h(\theta_0'x_1-\theta_0'x_2)K_h(y_1-y_2)f'^{\tau}_{\theta_0'X}(\theta_0'x_1)}{f^{\tau}_{\theta_0'X}(\theta_0'x_1)^2},
\end{align*}
where $f'^{\tau}_{\theta_0'X}$ denotes the derivative of $u\rightarrow f^{\tau}_{\theta_0'X}(u),$ the conditional density of $\theta_0'X$ given $Y\in A_{\tau}.$
\end{lemma}

\begin{lemma}
\label{leplusdur2} Under the Assumptions of Lemma \ref{parametric}
$$\sup_{h,\tau}\int\frac{[\nabla_{\theta}\hat{f}^{h,\tau}_{\theta_0}(y,x)-\nabla_{\theta}f^{*h,\tau}_{\theta_0}(y,x)]\mathds{1}_{y\in A_{\tau}}J_0(x,c/4)d\mathbb{P}(x,y)}{f_{\theta_0}^{{ \tau}}(y,\theta_0'x)}=o_P(n^{-1/2}).$$
\end{lemma}

\begin{proof}[Proof of Lemma \ref{leplusdur}]
To prove (\ref{01}-\ref{03}), we only prove (\ref{03}) since the others are similar. To prove (\ref{03}), we only consider the terms in which the second derivative is involved, the others being studied analogously.
Consider
\begin{align*}
&\frac{1}{h}\sum_{i=1}^n \delta_iW_{in}(X_i-x)K_h''(\theta'X_i-\theta'x)K_h(Z_i-y)(X_i-x)'\mathds{1}_{Z_i\in A_{\tau}}f^{\tau}_{\theta'X}(\theta'x)^{-1}\\
&=\frac{1}{h}\sum_{i=1}^n \delta_iW_{i}^*J_{\theta}(X_i,c)(X_i-x)K_h''(\theta'X_i-\theta'x)K_h(Z_i-y)(X_i-x)'\mathds{1}_{Z_i\in A_{\tau}}f^{\tau}_{\theta'X}(\theta'x)^{-1}\\
&+\frac{1}{h}\sum_{i=1}^n
\delta_iW_{i}^*Z_G(Z_i-)(X_i-x)K_h''(\theta'X_i-\theta'x)K_h(Z_i-y)(X_i-x)'\mathds{1}_{Z_i\in
A_{\tau}}f^{\tau}_{\theta'X}(\theta'x)^{-1},
\end{align*}
where the first term is contained in $\nabla^2_{\theta}f_{\theta}^{*h,\tau},$ while the second can be bounded by
$$O_P(n^{-1/2}h^{-2})\left[\frac{1}{nh^2}\sum_{i=1}^n\delta_i\mathds{1}_{Z_i\leq \tau_0}|K''|\left(\frac{\theta'X_i-\theta'x}{h}\right)|K|\left(\frac{Z_i-y}{h}\right)\right].$$ Using the results of \cite{Sherman94}, the term inside the brackets is $O_P(1)$ uniformly in $x,$ $y$, $\theta$ and $h.$

Now, for the representation (\ref{04}), observe that
\begin{align}
&\nabla_{\theta}[\hat{f}_{\theta_0}^{h,\tau}-f^{*h,\tau}_{\theta_0}](y,x)\nonumber\\
&\quad =h^{-1}\sum_{i=1}^n \delta_i W_i^*Z_G(Z_i-)(x-X_i)K_h'(\theta_0'x-\theta_0'X_i)K_h(y-Z_i)f^{\tau}_{\theta_0'X}(\theta_0'x)^{-1}\mathds{1}_{Z_i\in A_{\tau}} \nonumber \\
&\qquad-\sum_{i=1}^n \delta_i W_i^*Z_G(Z_i-)J_{0}(x,c)K(\theta_0'X_i-\theta_0'x)K_h(Z_i-y){f}'^{\tau}_{\theta_0'X}(\theta_0'x)f^{\tau}_{\theta_0'X}(\theta_0'x)^{-2}\mathds{1}_{Z_i\in A_{\tau}} \nonumber \\
&\qquad +R'_n(\tau,h,x,y), \label{decompopo}
\end{align}
with $\sup_{x,y,h,\tau}|R'_n(\tau,h,x,y)|=O_P\left(n^{-1}h^{-3/2}[\log n]^{1/2}\right),$ from the convergence rate of $Z_G$ (see (\ref{gill}) and (\ref{zhou})) and the convergence rate of the denominator in (\ref{kernel_estimator}) and its derivative, say $(\hat{f}^{\tau}_{\theta_0'X}-f^{\tau}_{\theta_0'X})$ and $(\hat{f}'^{\tau}_{\theta_0'X}-f'^{\tau}_{\theta_0'X})$ (which are of uniform rate $O_P\left(n^{-1/2}h^{-1/2}[\log n]^{1/2}\right)$ and $O_P\left(n^{-1/2}h^{-3/2}[\log n]^{1/2}\right)$ from arguments similar as for the proofs of (\ref{uniff})-(\ref{unif2}) and (\ref{01})-(\ref{03})). An i.i.d. representation of the main term in (\ref{decompopo}) can be deduced from Theorem \ref{km}
since the class $\{h^3g^{h,\tau}_{f,x,y},x,y,h\}$ is a VC-class from \cite{Nolan87}.
\end{proof}
\begin{proof}[Proof of Lemma \ref{leplusdur2}]
Observe that, from classical kernel arguments
\[\sup_{t}\left|\iint_{x_2,t \leq y_2 \leq \tau_0}g^{h,\tau}_{f,x,y}(x_2,y_2)J_0(x,c/4)d\mathbb{P}(x_2,y_2)d\mathbb{P}(x,y)-E[\nabla_{\theta}f_{\theta_0}^{{\tau}}(Y,X)J_0(X,c/4)]\right|=O(h^4),\] since $K$ is of order 4.
From the representation (\ref{04}) in Lemma \ref{leplusdur},
\begin{align}
&\int[\nabla_{\theta}\hat{f}^{h,\tau}_{\theta_0}(y,x)-\nabla_{\theta}f^{*h,\tau}_{\theta_0}(y,x)]J_0(x,c/4)d\mathbb{P}(x,y)\nonumber \\
&\quad=\int [1-H(t)]^{-1}E[\nabla_{\theta}f_{\theta_0}^{\tau}(Y,X)J_0(X,c/4)]d\bar{M}^G(t)\nonumber \\
&\qquad+\int [1-H(t)]^{-1}\left[\iint_{x_2,t \leq y_2 \leq \tau_0}g^{h,\tau}_{f,x,y}(x_2,y_2)J_0(x,c/4)d\mathbb{P}(x_2,y_2){d\mathbb{P}(x,y)}\right. \nonumber \\
&\qquad-E[\nabla_{\theta}f_{\theta_0}^{\tau}(Y,X)J_0(X,c/4)]\bigg]d\bar{M}^G(t)\nonumber \\
&\qquad+\int R_n(\tau,h,x,y)d\mathbb{P}(x,y), \label{encoreuneffort}
\end{align}
where the last term is $o_P(n^{-1/2})$ uniformly in $\tau$ and $h$. The first one is zero from Proposition \ref{expectation} and since $J_0$ only depends on $\theta_0'X.$

For the second, let $\phi_n(t,h,\tau)=$\\$[1-H(t)]^{-1}\{\iint_{x_2,t \leq y_2 \leq \tau_0}g^{h,\tau}_{f,x,y}(x_2,y_2)J_0(x,c/4)d\mathbb{P}(x_2,y_2)d\mathbb{P}(x,y)-E[\nabla_{\theta}f_{\theta_0}^{\tau}(Y,X)J_0(X,c/4)]\}.$
Using the fact that $\mathcal{H}_n$ is of cardinality $k_n,$ we have, for the second term in (\ref{encoreuneffort}),
\[
\mathbb{P}\left(\sup_{h\in \mathcal{H}_n}\left|\int \phi_n(t,h,\tau)d\bar{M}^G(t)\right|\geq \varepsilon\right) \leq  k_n \sup_{h\in \mathcal{H}_n}\mathbb{P}\left(\left|\int \phi_n(t,h,\tau)d\bar{M}^G(t)\right|\geq \varepsilon\right).
\]
Now apply Lenglart's inequality (see Lenglart (1977) or Theorem 3.4.1 in Fleming and Harrington, (1991)). This shows that, for all $\varepsilon>0$ and all $\eta>0,$
\begin{align}
\label{leng1}
\nonumber &\mathbb{P}\left(\sup_{\tau \leq s\leq \tau_0}\left\{\int_0^s \phi_n(t,h,\tau)d\bar{M}^G(t)\right\}^2\geq \varepsilon^2\right)\\
&\quad\leq \frac{\eta}{\varepsilon^2}+\mathbb{P}\left(n^{-1}\int_0^{\tau_0}\phi_n^2(t,h,\tau)\frac{[1-\hat{H}(t-)]dG(t)}{1-G(t-)}\geq \eta\right).
\end{align}
As mentioned before, $\sup_{t}|\phi_n(t,h,\tau)|=O(h^4).$ From (\ref{leng1}) and condition on $k_n$ in Assumption \ref{noyau}, the Lemma follows.
\end{proof}

\subsubsection{Donsker classes}

As stated in Assumption \ref{donsker1}, to obtain a $n^{-1/2}-$convergence of $\hat{\theta},$ we need the regression function (and its gradient) to be sufficiently regular. In the Lemma below, we first show that the classes of functions defined in Assumption \ref{donsker1} are Donsker, and that $\hat{f}^{h,\tau}_{\theta_0}$ also belongs to the same regular class as $f^{\tau}_{\theta_0}$ with probability tending to one.

\begin{proposition}
\label{classe de Donsker} Consider the classes $\mathcal{H}_1$ and $\mathcal{H}_2$ defined in Assumption \ref{donsker1}.
$\mathcal{H}_1$ and $\mathcal{H}_2$ are Donsker classes. Furthermore, $\hat{f}^{h,\tau}_{\theta_0}$
and $\nabla_{\theta}\hat{f}^{h,\tau}_{\theta_0}$ belong respectively to
$\mathcal{H}_1$ and $\mathcal{H}_2$ with probability tending to
one for some constant $M$ sufficiently large.
\end{proposition}

\begin{proof}
The class $\mathcal{H}_1$ is Donsker from Corollary 2.7.4 in \cite{Vaart96}. The
class $\mathcal{H}_2$ is Donsker from a permanence property of
Donsker classes, see Examples 2.10.10 and 2.10.7 in \cite{Vaart96}.
We only show the proof for
$\nabla_{\theta}\hat{f}^{h,\tau}_{\theta_0},$ since the one for
$\hat{f}^{h,\tau}_{\theta_0}$ is similar. Write
\begin{align*}
& \nabla_{\theta}\hat{f}^{h,\tau}_{\theta_0}(z,x) \\&=
\frac{1}{nh}\sum_{i=1}^n
\frac{\delta_i\mathds{1}_{Z_i \in A_{\tau}}(X_i-x)K'_h(\theta_0'X_i-\theta_0'x)K_h(Z_i-z)}{[1-\hat{G}(Z_i-)]f^{\tau}_{\theta_0'X}(\theta_0'x)}J_0(X_i,c/2)
\\
& \quad +\frac{1}{nh}\sum_{i=1}^n
\frac{\delta_i\mathds{1}_{Z_i \in A_{\tau}}(X_i-x)K'_h(\theta_0'X_i-\theta_0'x)K_h(Z_i-z)[\hat{f}^{\tau}_{\theta_0'X}(\theta_0'x)-f^{\tau}_{\theta_0'X}(\theta_0'x)]}{[1-\hat{G}(Z_i-)]\hat{f}^{\tau}_{\theta_0'X}(\theta_0'x)f^{\tau}_{\theta_0'X}(\theta_0'x)}J_0(X_i,c/2)
\\
& \quad -\left[\frac{1}{nh}\sum_{i=1}^n
\frac{(X_i-x)K_h'(\theta_0'X_i-\theta_0'x)J_0(X_i,c/2)}{\big(f_{\theta_0'X}^{\tau}(\theta_0'x)\big)^2}\right]
\\
& \quad \times \left[\frac{1}{n}\sum_{i=1}^n
\frac{\delta_iK_h(\theta_0'X_i-\theta_0'x)K_h(Z_i-z)\mathds{1}_{Z_i \in A_{\tau}}}{[1-\hat{G}(Z_i-)]}\right]
\\
& \quad +\left[\frac{1}{nh}\sum_{i=1}^n
\frac{(X_i-x)K_h'(\theta_0'X_i-\theta_0'x)\left[\big(\hat{f}_{\theta_0'X}^{\tau}(\theta_0'x)\big)^2-\big(f^{\tau}_{\theta_0'X}(\theta_0'x)\big)^2\right]J_0(X_i,c/2)}{\big(\hat{f}^{\tau}_{\theta_0'X}(\theta_0'x)f^{\tau}_{\theta_0'X}(\theta_0'x)\big)^2}\right]
\\ & \quad \times \left[\frac{1}{n}\sum_{i=1}^n
\frac{\delta_i\mathds{1}_{Z_i \in A_{\tau}}K_h(\theta_0'X_i-\theta_0'x)K_h(Z_i-z)}{[1-G(Z_i-)]}\right].
\end{align*}
From this expression, we clearly see that
$\nabla_{\theta}\hat{f}^{h,\tau}_{\theta_0}(y,x)=x\phi_1(x'\theta_0,y)+\phi_2(x'\theta_0,y).$
Now we must check that $\phi_1$ and $\phi_2$ are in
$\mathcal{H}_1$ with probability tending to one. Since the functions are
twice continuously differentiable (from the assumptions on $K$),
we only have to check their boundedness. From Lemma \ref{leplusdur}, this can be done at first by replacing $\hat{f}^{h,\tau}$ by $f^{*h,\tau}$ (i.e. $\hat{G}$ by the true function $G$). Among the several terms in the decomposition of $\nabla_{\theta}f^{*h,\tau},$ we will only study
$$\phi(u,y) = \frac{1}{nh}\sum_{i=1}^n
\frac{\delta_i\mathds{1}_{Z_i \in A_{\tau}}X_iK'_h(\theta_0'X_i-u)K_h(Z_i-z)J_0(X_i,c/2)}{[1-G(Z_i-)]f^{\tau}_{\theta_0'X}(u)},$$ since the others are similar.
We will show that the derivatives of order 0, 1 and $1+\delta$ of
this function are uniformly bounded by some constant $M$ with probability
tending to one.

Now a centered version of $\phi$ converges to zero at rate
$O_P([\log n]^{1/2}n^{-1/2}h^{-1})$ (see \cite{Einmahl05}), which tends to zero as long as
$nh^2\rightarrow \infty.$ Furthermore, $E[\phi]$ is uniformly bounded from our Assumption \ref{donsker1} on the regression function.
For the derivative,
\begin{align*}
\partial_u\phi(u,y) &=- \frac{1}{nh}\sum_{i=1}^n
\frac{\delta_i\mathds{1}_{Z_i \in A_{\tau}}X_iK''_h(\theta_0'X_i-u)K_h(Z_i-z)J_0(X_i,c/4)}{[1-G(Z_i-)]f^{\tau}_{\theta_0'X}(u)}\\
 & \quad -\frac{1}{nh}\sum_{i=1}^n
\frac{\delta_i\mathds{1}_{Z_i \in A_{\tau}}X_iK'_h(\theta_0'X_i-u)K_h(Z_i-z)J_0(X_i,c/4){f}'^{\tau}_{\theta_0'X}(u)}{[1-G(Z_i-)]\big(f^{\tau}_{\theta_0'X}(u)\big)^2}.
\end{align*}
Again, $E[\partial_u\phi]$ is uniformly bounded from our Assumption \ref{donsker1}.
Now using the results of \cite{Einmahl05}, the centered version of $\partial_u\phi$ tends to zero provided that
$nh^6\rightarrow \infty.$ The same arguments apply for $\partial_y\phi.$ Hence, with $f_i(u,y)=E(\phi_i(u,y))$ we proved that $\sup_{u,y} |\partial_u^{j}\partial_y^{k}\phi_{i}(u,y)-\partial_u^{j}\partial_y^{k}f_{i}(u,y)|$ tends to zero in probability for $i=1,2$, $k+j\leq 1$.
Now we have to show that $\partial_u\phi_j$ and $\partial_y\phi_j$
are $\delta-$ Hölder for $j=1,2$ with an H\"{o}lderian constant bounded by some $M$ with probability tending to one. We only prove the result for
$\partial_u\phi_1.$ We have
\begin{align*}
\sup_{u',y',x,y}\frac{\left|\partial_u\phi_{1}(u,y)-\partial_u\phi_{1}(u',y')\right|}{\|(u,y)-(u',y')\|^{\delta}}
&= \max\left(\sup_{|u-u'|\geq n^{-1}
,y,y'}\frac{\left|\partial_u\phi_{1}(u,y)-\partial_u\phi_{2}(u',y')\right|}{\|(u,y)-(u',y')\|^{\delta}},\right.
\\
& \left.\sup_{|u-u'|\leq n^{-1}
,y,y'}\frac{\left|\partial_u\phi_{1}(u,y)-\partial_u\phi_{1}(u',y')\right|}{\|(u,y)-(u',y')\|^{\delta}}\right)
\\
&= \max(S_1,S_2).
\end{align*}
We have
\begin{align*}
S_1 &\leq
\sup_{u,y,u',y'}\frac{|\partial_uf_{1}(u',y')-\partial_uf_1(u,y)|}{\|(u',y')-(u,y)\|^{\delta}}
\\
&
+2n^{\delta}\sup_{u,y,u',y'}|\partial_u\phi_{1}(u,y)-\partial_uf_1(u,y)|.
\end{align*}
From our Assumptions, the first supremum is bounded, while the
last is\\
$O_P(n^{-1/2+\delta}[\log n]^{1/2}h^{-3})$ from the
convergence rate of $\partial_u\phi_2.$ It tends to zero provided
that
$nh^{6+\delta}\rightarrow \infty.$
For $S_2,$ since $K$ is $\mathcal{C}^3$ with bounded derivatives, for some positive constant
$M,$
\begin{align*}
\sup_{\|(u,y)-(u',y')\|\leq n^{-1}
,y,y'}\frac{\left|\partial_u\phi_{1}(u,y)-\partial_u\phi_{1}(u',y')\right|}{\|(u,y)-(u',y')\|^{\delta}}
\leq M\times
O_P(1)\|\sum_{i=1}^3|K^{(i)}|\|_{\infty}\\
\times \sup_{\|(u,y)-(u',y')\|\leq n^{-1}
}\|(u,y)-(u',y')\|^{1-\delta}h^{-1}\frac{1}{nh^4}\sum_{i=1}^n\frac{\delta_i}{1-G(Z_i-)}.
\end{align*}
The last supremum is bounded by $O_P(1)\times
n^{-1+\delta}h^{-5},$ and it tends to zero when $nh^6\rightarrow
\infty$ (and the $O_P(1)$ term does not depend on $u,y$).
\end{proof}

\end{document}